\documentclass[onecolumn,11pt]{article}
\usepackage[top=1in, bottom=1in, left=1in, right=1in]{geometry}
\usepackage{amsmath,amsfonts,amscd,amssymb}
\usepackage{graphicx}
\usepackage{epstopdf}
\usepackage{overpic}
\usepackage{cancel}
\usepackage{rotating}
\usepackage{url}
\usepackage{caption}
\usepackage{color}
\usepackage{rotating}
\usepackage{multirow}
\usepackage{wrapfig}
\usepackage{mathtools}
\usepackage{subeqnarray}
\usepackage{setspace}
\usepackage{palatino} 
\usepackage[numbers,sort&compress]{natbib}

\usepackage[bottom,flushmargin,hang,multiple]{footmisc}
\usepackage{lipsum}
\usepackage{chngcntr}

\newcommand\blfootnote[1]{%
  \begingroup
  \renewcommand\thefootnote{}\footnote{#1}%
  \addtocounter{footnote}{-1}%
  \endgroup
}

\definecolor{header1}{cmyk}{0,0,0,1}

\newcommand{\bv}{\mathbf{v}}
\newcommand{\bx}{\mathbf{x}}

\newcommand{\bu}{\mathbf{u}}

\newcommand{\bg}{\mathbf{g}}

\title{\Large{\vspace{-.55in}\textbf{Finite Time Lyapunov Exponent Analysis\\ of Model Predictive Control and Reinforcement Learning}}\vspace{-.175in}}

\author{\normalsize{Kartik Krishna$^{1*}$, Steven L. Brunton$^1$, and Zhuoyuan Song$^2$}\\
\footnotesize{$^1$ Department of Mechanical Engineering, University of Washington, Seattle, WA 98195, United States}\\
\footnotesize{$^2$ Department of Mechanical Engineering, University of Hawai`i at M\={a}noa, Honolulu, HI 96822, United States}\vspace{-.1in}}
\date{}
\begin{document}
\maketitle

\blfootnote{$^*$ Corresponding author (karkris3@uw.edu).}
\vspace{-.2in}
\begin{abstract}
 Finite-time Lyapunov exponents (FTLEs) provide a powerful approach to compute time-varying analogs of invariant manifolds in unsteady fluid flow fields.
These manifolds are useful to visualize the transport mechanisms of passive tracers advecting with the flow.
However, many vehicles and mobile sensors are not passive, but are instead actuated according to some intelligent trajectory planning or control law; for example, model predictive control and reinforcement learning are often used to design energy-efficient trajectories in a dynamically changing background flow.
In this work, we investigate the use of FTLE on such controlled agents to gain insight into optimal transport routes for navigation in known unsteady flows.
We find that these controlled FTLE (cFTLE) coherent structures separate the flow field into different regions with similar costs of transport to the goal location.
These separatrices are functions of the planning algorithm's hyper-parameters, such as the optimization time horizon and the cost of actuation. 
Computing the invariant sets and manifolds of active agent dynamics in dynamic flow fields is useful in the context of robust motion control, hyperparameter tuning, and determining safe and collision-free trajectories for autonomous systems.
Moreover, these cFTLE structures provide insight into effective deployment locations for mobile agents with actuation and energy constraints to traverse the ocean or atmosphere.
\\

\noindent\emph{Keywords-optimal control, finite time Lyapunov exponents, path planning, mobile sensors, dynamical systems, unsteady fluid dynamics, model predictive control, reinforcement learning}
\end{abstract}

\section{Introduction}
Trajectory planning in an unsteady flow field is an important problem for intelligent mobile agents, with applications including environmental monitoring and data collection~\cite{fossum2019toward,chai2020monitoring,zhang2021system,bellingham2007robotics,wynn2014autonomous,rhoads2013minimum}. 
When planning trajectories, many applications aim at achieving certain objectives ranging from reaching a static goal location to maintaining certain connectivity of a multi-agent sensor network for part of or the entire the mission~\cite{song2017multi, song2015anisotropic}.
Optimization and control are often employed in designing the decision-making algorithms on-board the mobile agents, enabling offline or real-time trajectory planning to achieve the desired objectives.
Intelligent algorithms that leverage the background flow are necessary, since naively using full propulsion while aiming at a target can result in wasteful trajectories and the potential of the vehicle being swept away by large currents at a later time.
However, even with on-board algorithms, it is still imperative to carefully choose the deployment locations since the agent's ability to reach certain regions is largely determined by its actuation limits and the background flow dynamics.
For example, it might be impossible for two groups of agents that are dominated by close-by, but different flow structures, to rendezvous.
Furthermore, tuning the hyperparameters of an on-board control strategy to obtain the best performance is a challenging task. 
The ability to summarize and visualize the dependence of the control performance on the control hyperparameters may aid in this process.
In this work, we investigate the use of finite-time Lyapunov exponents (FTLEs) from dynamical systems to address these challenges by quantifying the performance and sensitivity of planning and control algorithms.
We also discover a mathematical connection between optimal control and FTLEs of the controlled flow.

In dynamical systems, Lyapunov exponents provide a measure of the sensitivity of the trajectory to initial conditions. 
In chaotic vector fields, two almost identical initial conditions can lead to flows that diverge exponentially in finite time. 
The finite-time Lyapunov exponent quantifies this stretching over a fixed, finite-time horizon, resulting in a scalar field over a domain of interest highlighting the most sensitive regions to perturbations in the initial conditions~\cite{shadden2005physd,Green2007jfm, brunton2010chaos,lipinski:2010,haller2015arfm}. 
Moreover, the FTLE has also been used in fluid dynamics applications to compute Lagrangian coherent structures (LCS), which are finite-time analogues of invariant manifolds that mediate the transport of material in unsteady fluid flows~\cite{haller2002pof,haller:05,shadden2005physd,Green2007jfm,shadden2009correlation,shadden2011lagrangian,haller2015arfm,sudharsan2016lagrangian}. 
The LCS, and consequently the FTLE, define transport barriers in a flow field where passive drifters are attracted to or repelled from.
The FTLE method has been successfully applied to bio-propulsion~\cite{wilson2009lagrangian}, medicine~\cite{shadden2008characterization,forgoston2011maximal}, the spread of microbes~\cite{tallapragada2011lagrangian}, and the study of aerodynamics~\cite{rockwood2019practical,rockwood2019real}, among other domains.

For mobile agents with actuation capabilities, there exists a wealth of knowledge on various algorithms for trajectory generation in dynamic fluid environments. 
For example, graph search algorithms and stochastic optimization have been investigated for path planning~\cite{kularatne2016time, rao2009large, subramani2016energy}. 
Assimilating in-situ observations obtained by mobile sensors in an adaptive fashion into ocean models has also been explored, for example with mixed integer programming algorithms~\cite{yilmaz2008path, lermusiaux2007adaptive}. 
Coordinated control of ocean gliders for adaptive ocean sensing has been exhaustively studied~\cite{bhatta2005coordination, leonard2007collective, fiorelli2006multi, leonard2001model}. 
Algorithms inspired from computational fluid dynamics have also been used to explore coordinated control of swarms in flow fields~\cite{lipinski2010cooperative,lipinski2014feasible, lipinski2011master, song2017multi, song2015anisotropic}. 
Recent developments have made use of model predictive control (MPC)~\cite{krishna2022finite} and reinforcement learning (RL) algorithms~\cite{gunnarson2021learning, bifarale2019, buzzicotti2021optimal, jiao2021learning,ao2023} to find optimal paths in unsteady flow fields.
Optimal control, including MPC, has been previously related to the passive FTLE in the past~\cite{inanc2005optimal, zhang2008optimal,senatore2008fuel,heckman2016controlling}. 
Among the existing methods, model predictive control~\cite{garcia1989model,camacho2013model} and reinforcement learning~\cite{sutton2018reinforcement,recht2019tour} are among the most useful paradigms in modern control~\cite{Brunton2022book}, and analyzing these control laws and \emph{policies} with FTLE be the focus of this paper.

A control policy, for example from MPC or RL, maps each agent state in the domain of fluid flow to an action, which can ultimately be visualized as a vector field. 
When agents follow a policy within an unsteady flow field, their behaviour can be understood as passive drifters operating within an entirely new active flow field -- one that is a combination of the original background flow field and the control policy.
Given the applicability of FTLE for understanding \emph{passive}, uncontrolled transport in unsteady fluid flow fields, it is natural to extend it to understand \emph{active} transport given an agent's control policy, which is the focus of this paper. 
Policies derived from optimal control methods, and by extension the resulting active flow fields, can often be spatially discontinuous or non-smooth depending on the hyperparameters used.
Therefore, our work also contributes to understanding the use of FTLE on non-smooth systems.

\begin{figure}[t]
    \centerline{\includegraphics[width = \linewidth,trim = {0 150mm 0 60mm }, clip, ]{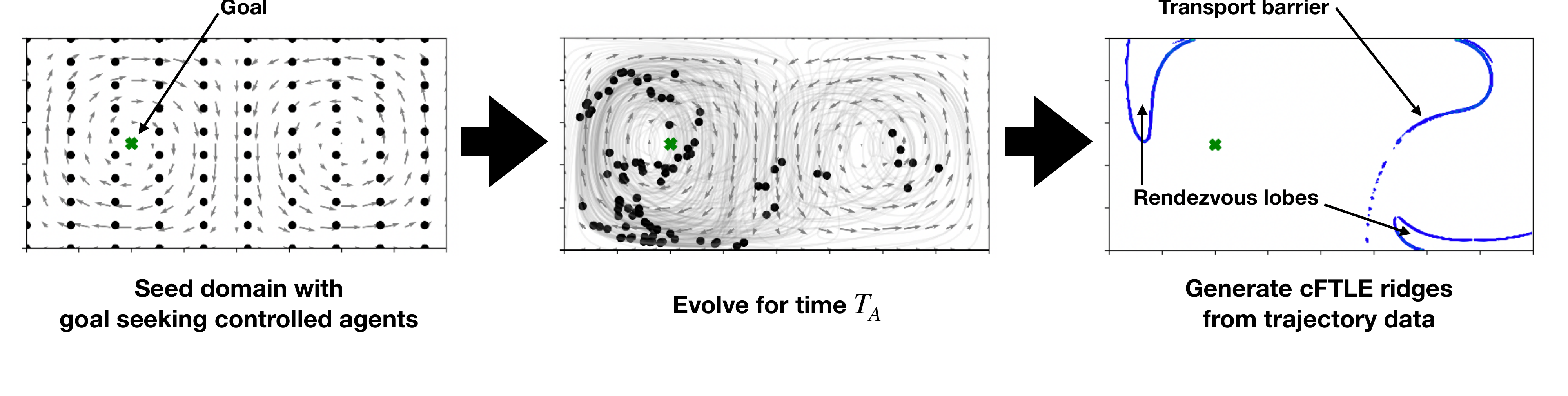}}
    \vspace{-.25in}
    \caption{This figure outlines the methodology of cFTLE which takes Lagrangian controlled trajectories as input and outputs LCS curves which highlight interesting transport features of the control law which are generally not visualizable by simply plotting all the trajectories or by plotting the control vector field. The left most figure shows that in order to evaluate our control algorithm, we generate the simulation from a mesh grid of initial agent positions at $t=0$. We end the simulations at  $t=T_A$ (as shows in the middle figure). Finally we use the FTLE algorithm on this data to generate “cFTLE” ridges which highlight important features in the domain. }
    \label{fig:overview}
\end{figure}

In this paper, we compute FTLE on active agents navigating within an unsteady fluid flow field. 
Specifically, we examine the controlled trajectories using either a finite-horizon model predictive control optimization or a reinforcement learning strategy to learn a policy that drives an agent towards a goal state.
We then perform control FTLE (cFTLE) analysis on the combined flow field generated from the background flow and the control policy, illustrated in Fig.~\ref{fig:overview}.
Similar studies have been performed in the past for biological applications, where FTLE was used to understand the behavior of plankton actively propelling away from feeding jellyfish~\cite{peng2009transport}.
We confirm that, much like the traditional FTLE, cFTLE can be used to visualize transport barriers, which allow us to understand the transport of agents following a control policy.
These results are useful in the context of finding ideal deployment locations for rendezvous and coordinated multi-agent goal tracking.
We further relate the features of cFTLE, which uncover exponential divergence of initial conditions, to the sensitivity of optimization cost functions used in trajectory planning. We also show that these findings are not exclusive to policies generated through MPC or RL, but to policies generated by any optimal control method with a particular structure in cost function, which we elaborate on.
In a previous work by the authors~\cite{krishna2022finite}, it was found that large energy expenditure of active agents corresponded with the presence of large background FTLE.
Here, we further explore this connection and find that the energy spent has a closer correlation with cFTLE, which are deformations of the passive FTLE ridges.
Given the connections to optimal control, we perform sweep through a range of parameters that govern the aggressiveness of the control, the set point to track, and the amount of future information about the background flow field available, which creates spatial discontinuities. 
In the future, this study could potentially aid in methods for faster learning of control policies with partial knowledge of the control FTLE and vice versa, and aid in summarizing the effectiveness of policies used for flow navigation.

\section{Background}
The main idea of this paper is the use of the finite-time Lyapunov exponents to analyze active agents moving towards a goal within an unsteady fluid flow field. 
Traditionally, FTLE analysis has been primarily used on trajectories of passive drifters to highlight coherent structures within the flow field that mediate transport of the drifters.
In this section, we give a short background on the classical passive FTLE method.
MPC is used extensively in this work to establish connections between energy-optimal trajectories and the FTLE. 
However, these connections are not unique to MPC; policies generated through other formulations of optimal control can also be related to the FTLE, as we will show with the case of reinforcement learning. 
With this in mind, we broadly introduce the reader to optimal control theory and notions of the value function in this section.

\subsection{Finite Time Lyapunov Exponents}
The original FTLE~\cite{haller2002pof,shadden2005physd,haller2015arfm} primarily addressed dynamical systems of the form
\begin{equation}
 \frac{d }{dt}\bx(t)  = \mathbf{f} \left(\bx(t),t\right),
 \label{eq:uncontrolled}
\end{equation}
where the function $\mathbf{f} \left(\bx(t),t\right):\mathbb{R}^n \times \mathbb{R} \to  \mathbb{R}^n$ represents the system dynamics, $t \in \mathbb{R}$ represents time, and $\bx(t) \in \mathbb{R}^n$ is the state of the system.
An FTLE algorithm takes a dynamical system as input and generates a scalar variable that can be used in the computation of separatrices, known as Lagrangian coherent structures (LCS) in the flow field. These separatrices demarcate the boundaries between different regions in a domain of interest where passive tracers remain trapped~\cite{kelley2013lagrangian}. This shall be seen later in the paper, particularly highlighted in Figure~\ref{fig:regimes}. In the past, FTLE analysis has mostly been performed using passive drifters to study time-varying vector fields, such as ocean flows~\cite{olascoaga2006persistent,beron2008oceanic,beron2015dissipative} and pollution transport models~\cite{Lekien2005physicad}. 
More broadly, FTLE has also been used to compute coherent structures for a wide range of other flows~\cite{Green2007jfm,franco:2007,Padberg2007njp,Mathur2007prl,Peng2008jeb,rockwood2016detecting}. 

The FTLE field for the dynamical system $\mathbf{f} \left(\bx(t),t\right)$ can be computed as follows.
First, we initialize a grid of passive drifter particles at time $t_0$ and numerically integrate them through $\mathbf{f}(\bx(t),t)$ for a fixed amount of time (i.e., the advection time) $T_A \in \mathbb{R}$, resulting in a flow map ${\boldsymbol{\Phi}}^{t_0 + T_A}_{t_0}: \mathbb{R}^n \to  \mathbb{R}^n$:
\begin{equation}
    {\boldsymbol{\Phi}}^{t_0 + T_A}_{t_0}: \mathbf{x}(t_0) \mapsto \mathbf{x}(t_0) + \int_{t_0}^{t_0 +T_A} \mathbf{f}\left( \mathbf{x}(\tau),\tau\right) \, d\tau.
    \label{eq:integrate}
\end{equation}
 The operator ${\boldsymbol{\Phi}}^{t_0 + T_A}_{t_0}$ maps an initial state $\bx(t_0)$ to a final state $\bx(t_0 + T_A)$ by the differential equation flow induced in phase space.

Next, a Jacobian matrix of partial derivatives of the flow map, $\mathbf{D} {\boldsymbol{\Phi}}^{t_0 + T_A}_{t_0}$, is computed using finite differences for each drifter on the grid, represented by the grid node indices $i, j \in \mathbb{Z}^+$, such that 

\begin{equation}
\begin{aligned}
\left({\mathbf{D} {\boldsymbol{\Phi}}^{t_0 + T_A}_{t_0}}\right)_{i,j}
&{ =
\begin{bmatrix} 
\frac{\Delta x_i(t_0+T_A)}{\Delta x_i(t_0)} & \frac{\Delta x_j(t_0+T_A)}{\Delta y_j(t_0)} \\ \frac{\Delta y_i(t_0+T_A)}{\Delta x_i(t_0)} & \frac{\Delta y_j(t_0+T_A)}{\Delta y_j(t_0)}
\end{bmatrix}} \\
&=
\begin{bmatrix}
 \frac{x_{i+1,j}(t_0 + T_A) - x_{i-1,j}(t_0 + T_A)}{x_{i+1,j}(t_0) - x_{i-1,j}(t_0)} 
 & \frac{x_{i,j+1}(t_0 + T_A) - x_{i,j-1}(t_0 + T_A)}{y_{i,j+1}(t_0) - y_{i,j-1}(t_0)} \\
\frac{y_{i+1,j}(t_0 + T_A) - y_{i-1,j}(t_0 + T_A)}{x_{i+1,j}(t_0) - x_{i-1,j}(t_0)} 
& \frac{y_{i,j+1}(t_0 + T_A) - y_{i,j-1}(t_0 + T_A)}{y_{i,j+1}(t_0) - y_{i,j-1}(t_0)}
\end{bmatrix},
\end{aligned}
\end{equation}
where $x,y \in \mathbb{R}$ are the horizontal and vertical components of the position vector $\bx(t)$. The key idea here is that nearby initial conditions that rapidly separate from each other in finite advection time are highlighted as regions of large FTLE, representing repelling coherent structures. 
This flow map Jacobian is used to compute the Cauchy-Green deformation tensor
\begin{equation}
    \boldsymbol{\Delta}_{i,j} = \left({{\mathbf{D} {{\boldsymbol{\Phi}}^{t_0 + T_A}_{t_0}}}}\right)^* \mathbf{D}{{\boldsymbol{\Phi}}^{t_0 + T_A}_{t_0}},
\end{equation}
where $^*$ represents the matrix transpose. 
Finally, the largest eigenvalue $\lambda_{\text{max}}$ of $\boldsymbol{\Delta}_{i,j}$ for each drifter $i,j$ is used to compute the FTLE field:
\begin{equation}
   \sigma_{i,j} = \frac{1}{|T_A|} \ln {\sqrt{{(\lambda_{\text{max}})}_{i,j}}}.
\end{equation}
Alternatively, $\sigma_{i,j}$ can be viewed as the maximum singular value from the singular value decomposition (SVD) of $\mathbf{D}{\boldsymbol{\Phi}}^{t_0 + T_A}_{t_0}$. 
It is important to note that for unsteady flow fields, the FTLE field will also vary in time, so that at each new time step a new grid of drifters must be reinitialized and advected through the flow. 
This procedure is typically quite expensive to compute, although there are algorithms to speed up the calculations~\cite{brunton2010chaos,lipinski:2010}.

The ridges of the computed FTLE field can be extracted to visualize manifolds in the domain, known as Lagrangian coherent structures. This requires an additional step of computing the Hessian of $\sigma_{i,j}$ for ridge extraction.  
FTLE based on drifter particles integrated forward in time ($T_A>0$) result in coherent structures that repel particles. This can be seen as the blue curves in the right most plot of Figure~\ref{fig:overview}. 
Similarly, FTLE based on particles integrated backward in time ($T_A<0$) results in attracting coherent structures.  

\subsection{Optimal Control}
Equation~\eqref{eq:uncontrolled} may be modified to include actuation and control, resulting in the state-space equation
\begin{equation}
 \frac{d }{dt}\bx(t)  = \mathbf{g} \left(\bx(t),\bu,t\right),
 \label{eq:controlled}
\end{equation}
where $\mathbf{g} \left(\bx(t), \bu,t\right):\mathbb{R}^n \times \mathbb{R}^m \times \mathbb{R} \to  \mathbb{R}^n$ , $t \in \mathbb{R}$ represents time, and $\bx(t) \in \mathbb{R}^n$ is the state of the agent. The agent is also able to apply actuation (such as propulsion) that is modeled as $\bu \in \mathbb{R}^m$. 

In optimal control, one seeks to find a control policy $\mathbf{u}(\mathbf{x},t)$ that minimizes a cost function $J$. Intuitively, this function often penalizes two terms, the distance from the goal we wish to drive the system states toward and the energy spent. Therefore, when the control is found which minimizes this function, the agent moves toward the goal in an energy efficient manner. In this paper, we primarily consider a linear–quadratic regulator (LQR) type quadratic cost function given by
\begin{equation}
    J_{LQR}(\mathbf{x}_0,t) = \int_{t_0}^{t_0 + T_H} \left[\mathbf{e(\bx_0, \tau)}^T\mathbf{Q} \; \mathbf{e}(\bx_0,\tau) +  \mathbf{u(\tau)}^T\mathbf{R}\;\mathbf{u}(\tau) \right]d\tau,
    \label{eq:cost_fn1}    
\end{equation}
where $\bx_0 \in \mathbb{R}^n$ is the initial spatial location of an agent, $t_0 \in \mathbb{R}$ is the initial time,  $\mathbf{e}(\bx_0, t) \triangleq\bx (\bx_0, t) -\bx_{\text{goal}}$ is the distance from goal or the state tracking error starting from an initial condition $\bx_0$,  $\mathbf{u} \in \mathbb{R}^m$ is the control vector, $T_H \in \mathbb{R}$ is the time horizon, $\mathbf{Q} \in \mathbb{R}^{n \times n}$ is the state penalty matrix, and $\mathbf{R} \in \mathbb{R}^{m \times m}$ is a matrix of control penalty. 

Functions of the form described by \eqref{eq:cost_fn1} are \emph{value functions}~\cite{sutton2018reinforcement}, which can also be viewed theoretically as solutions to the Hamilton-Jacobi Bellman (HJB) equations~\cite{Brunton2022book}.
These functions assign a scalar value to each initial condition in the domain corresponding to the future cost that will be accrued from that particular initial condition following the optimal policy.
The field of deep reinforcement learning often involves estimating these functions using neural networks from limited data of the agent moving through the domain or from having limited-to-no knowledge of the governing dynamics.
Analytically, some useful mathematical relations have been established~\cite{stengel1994optimal, todorov2006optimal} when using cost functions of the above form, coupled with control-affine or kinematic models highlighted in Eq.~\eqref{eq:kinemat}. One major result is 
\begin{equation}
    \bu^*(\mathbf{x}_0,t) = -\mathbf{R}^{-1} \nabla  J_{LQR}(\mathbf{x}_0,t),
    \label{eq:ctrl_grad}    
\end{equation}
which directly relates the spatial gradient of the cost function in the domain to the \emph{optimal} control law $\bu^*(\mathbf{x}_0,t)$ that minimizes Eq.~\eqref{eq:cost_fn1}. We will use this result later to establish connections to FTLE. When $\bu(\bx,t)$ is plotted spatially, it can be viewed as a \emph{policy}.

\section{Methodology for Control FTLE}
In this paper, we propose the use of FTLE in the analysis of policies generated through methods such as MPC and reinforcement learning. 
However, generating policies in this manner presents obstacles in the computation of cFTLE ridges, which will be elaborated in the following subsections. 
In this section, we outline how the computation of cFTLE can be derived from the traditional direct application of FTLE on unsteady flow fields. We also discuss the models used, the parameter ranges considered, and the tools used for policy generation.

\subsection{Control Finite Time Lyapunov Exponents (cFTLE)}\label{sec:cftle}
Extending the FTLE from the previous section to incorporate control, the equations governing the dynamics of an actuated agent from Eq.~\eqref{eq:controlled} are used. 
In this paper, we shall consider the class of kinematic models~\cite{inanc2005optimal} where
\begin{equation}
    \mathbf{g} \left(\bx(t), \bu,t\right) = \bv \left(\bx(t),t\right) + \bu.
    \label{eq:kinemat}
\end{equation}
These models assume that the sensor can generate its own flow relative velocity $\bu(t)=[u_x, u_y] \in \mathbb{R}^2$ in addition to the flow-induced velocity. 
Alternatively, models that incorporate inertial effects can also be used~\cite{sudharsan2016lagrangian, zhang2008optimal, peng2009transport}.
The function $\bv\left(\bx(t),t\right)$ is the unsteady flow field within which the agent is moving.
For the majority of this study, we consider this function to be given by the non-autonomous double-gyre equations~\cite{JiangS:95a,SpeichS:94a,SpeichS:95a}. 
Because the control $\bu$ is typically a full-state feedback control law, $\bu(\bx,t)$, the resulting flow field $\bg$ may be considered to be a function of $\bx$ for a given policy $\bu$.

Similar to traditional FTLE, we initialize a grid of agents at time $t_0$ and numerically integrate through $\mathbf{g}(\bx(t), \bu,t)$ for a fixed amount of time (i.e., the advection time) $T_A \in \mathbb{R}$, resulting in a flow map ${\boldsymbol{\hat\Phi}}^{t_0 + T_A}_{t_0}: \mathbb{R}^n \times \mathbb{R}^m \to  \mathbb{R}^n$:
\begin{equation}
    {\boldsymbol{\hat\Phi}}^{t_0 + T_A}_{t_0}: \mathbf{x}(t_0) \mapsto \mathbf{x}(t_0) + \int_{t_0}^{t_0 +T_A} \mathbf{g}\left( \mathbf{x}(\tau),\bu,\tau\right) \, d\tau.
    \label{eq:integrate}
\end{equation}
The flow map operator is particularly important in this paper as it will later be used to understand the connections between cFTLE and value functions in optimal control. The operator ${\boldsymbol{\hat\Phi}}^{t_0 + T_A}_{t_0}$ maps an initial agent position $\bx(t_0)$ to a position advected forward by time to $\bx(t_0 + T_A)$ by the flow field and agent propulsion.

FTLE based on agents integrated forward in time, $T_A>0$, result in coherent structures that repel agents, which can be seen in the right most plot in Figure~\ref{fig:overview} where the blue curves are visible, and later figures in this paper. 
Compared to the backward time FTLE, the forward time cFTLE is more straightforward, assuming that the policy for all states is known beforehand. 
To calculate the backward time cFTLE, we need both flow field and policy data backwards in time, which is not possible in realistic scenarios.
FTLE based on agents integrated backward in time, $T_A<0$, results in coherent structures which attract agents. 
These can be computed for periodic flow fields as both policy and flow field can be extrapolated backwards due to periodicity (as demonstrated in the Appendix).

\subsection{Problem Setup}
By combining the mobile sensor model and the double gyre flow field, the dynamics of the sensor are given by, 
\begin{equation}
\frac{d}{dt} 
\begin{bmatrix}
x\\y
\end{bmatrix}
= 
\begin{bmatrix}
-\pi A \sin(\pi f(x,t))\cos(\pi y)  \\
\pi A\cos(\pi f(x,t))\sin(\pi y)
\end{bmatrix}
+ 
\begin{bmatrix}
u_x\\u_y
\end{bmatrix}.
\label{eq:final_model}
\end{equation}
The time dependency is introduced by
\begin{align}
f(x,t) = a(t)x^2+b(t)x,
\end{align}
with time dependent coefficients
\begin{align}
a(t) = \epsilon \sin (\omega t) \quad \text{and} \quad
b(t) = 1-2\epsilon \sin(\omega t). \nonumber
\end{align}
This flow is defined on a nondimensionalized domain of $[0,2] \times [0,1]$, where $A, \epsilon, \omega \in \mathbb{R}$ are model parameters that control the flow pattern. 
Here, $\epsilon$ dictates the magnitude of oscillation in the $x$-direction, $\omega$ is the angular oscillation frequency, and $A$ controls the velocity magnitude. 
Unless stated otherwise, the parameters used for the double gyre flow field are as in Shadden et al.~\cite{shadden2005physd}, where $A = 0.1$, $\epsilon = 0.25$, and $\omega = 2 \pi / 10$.

For each agent, the control algorithm computes $\bu(\bx,t)$ numerically for a grid of initial conditions in the domain. 
{To solve the resulting optimization problems of MPC}, we use the CasADi~\cite{Andersson2019} and MPCTools~\cite{rawlings2015} packages, {which use an interior point filter-line search algorithm (IPOPT). }
{The discretized time step used for RK4 methods and for direct multiple shooting collocation in MPC across all double gyre simulations is $\Delta t = 0.1$.}
In later advanced examples with reinforcement learning, we use the stable baselines package~\cite{stable-baselines3}.

A relevant aspect of MPC is that it computes the optimization over a finite time horizon $T_H$, which is similar to how FTLE performs its computations over a time of advection $T_A$.
The objective in most examples discussed in the paper is to move the sensor from a starting location to a goal location at $\bx_{\text{goal}} = [0.5, 0.5]$. 
The optimal control problem is solved numerically with added constraints,
\begin{equation*}
    |u_x| \leq 0.1 \quad \text{and} \quad |u_y| \leq 0.1,
\end{equation*}
which ensure that the maximum sensor velocity is significantly smaller than the largest background flow field velocity, $\pi A\approx 0.314$. 
This constraint is imposed to model the limited actuation available in real world scenarios. 
From Eq.~\eqref{eq:cost_fn1}, in our case, $\mathbf Q = Q\mathbf{I}_{2 \times 2}$ and $\mathbf{R} = R\mathbf{I}_{2 \times 2}$, where $\mathbf{I}$ is the identity matrix. 
We combine these two parameters into one parameter, namely the $R/Q$, the specific values for which will be elaborated on as they appear in the text.

\subsection{Computation Using Interpolation}
For the computation of cFTLE, Eqs.~\eqref{eq:controlled}~and~\eqref{eq:integrate} often require the use of a dense grid of initial conditions.
Generally, with MPC, this would incur a large computational cost as a full \emph{closed-loop} MPC optimization would need to be computed over the prescribed time horizon, $T_H$, for the advection time, $T_A$, for each spatial initial condition at each simulation time step.
To reduce the total number of optimizations computed and reduce the requirement for several closed-loop MPC simulations, we pre-compute an open-loop policy for the entire space and time domain.
We consider agent initial conditions fixed to a grid (as in the left plot in Figure~\ref{fig:overview}) and compute the control actions at the grid node locations (as seen later in the plots of Figure~\ref{fig:policies}).
To keep the computational cost down, we make use of a coarse grid and then use linear interpolation to compute the policy for off-grid initial conditions.
This allows us to build a lookup-table function $\hat{\bu}(\bx, t)$ for every point in space and time.
We then integrate a dense grid of passive particles through the combined vector field of fluid flow and estimated \emph{open-loop } control through the following ODE,
\begin{equation}
 \frac{d }{dt}\bx(t)  = \mathbf{v} \left(\bx(t),t\right) + \hat{\bu},
 \label{eq:controlled_intp}
\end{equation}
where $\hat{\bu}$ is the estimated open loop control law using liner interpolation. 
We then follow the computation method outlined in Section~\ref{sec:cftle} to obtain the Lagrangian trajectories.

\begin{figure}[t]
    \centerline{\includegraphics[scale=.3,trim = {0 60mm 0 40mm }, clip,]{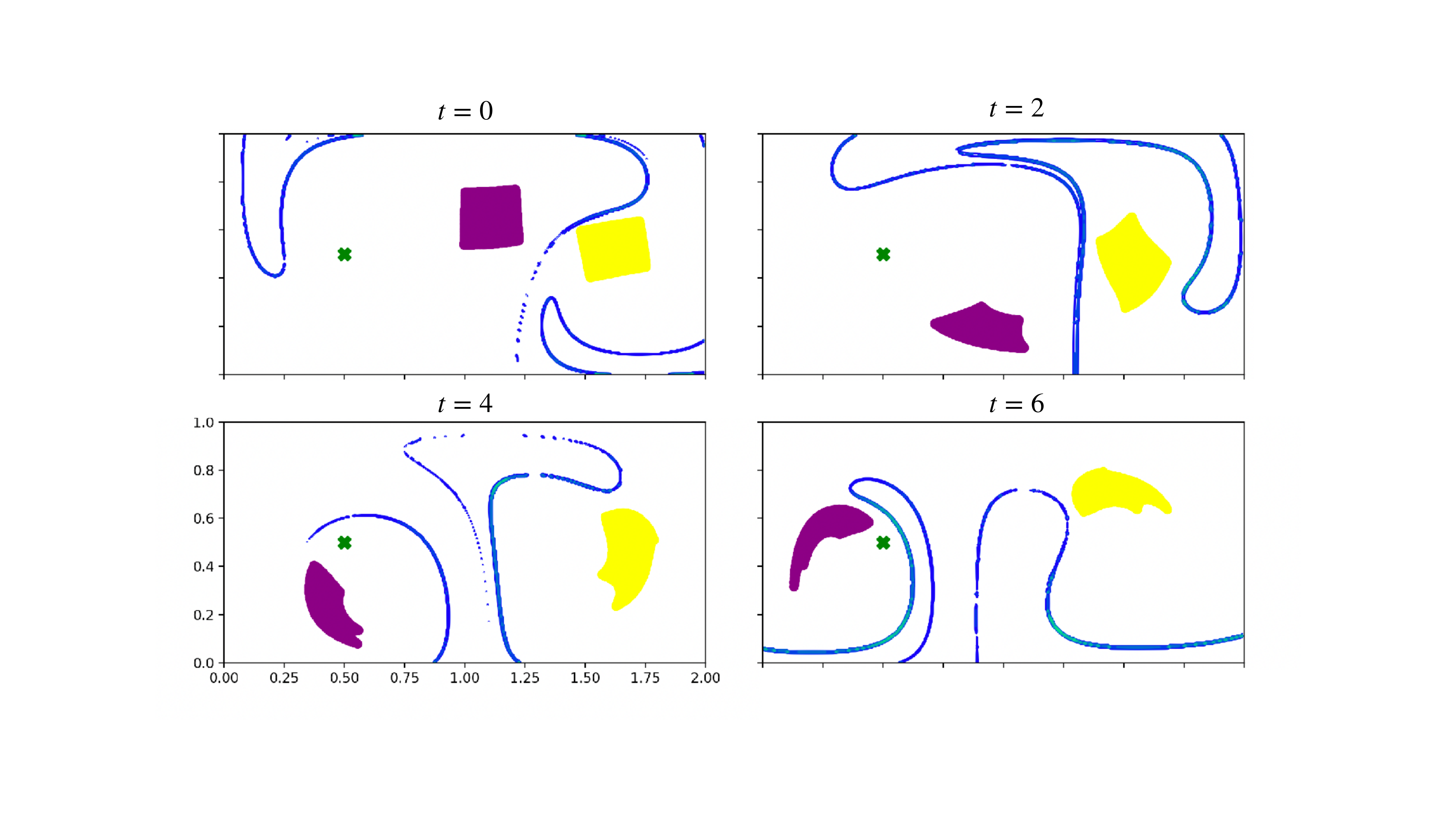}}
    \caption{Control FTLE ridges separate dynamically different regions of agents moving in the flow field. We show two distinct patches of agents - a purple and a yellow patch. Both patches start initially on different sides of the cFTLE ridge. The time evolution shows that cFTLE ridges separate the patches that undergo different paths and accrue largely different short-term state error cost despite starting close to each other.}
    \label{fig:regimes}
\end{figure}
 
\section{Results}
In this section, we first verify that cFTLE behaves much like regular FTLE and demonstrate how this property can be used diagnostically for the analysis of control laws. 
Next, we highlight a connection between value functions from optimal control theory and cFTLE ridges.
Finally, given the aforementioned connection, we study the deformation of cFTLE ridges under the change in parameters of MPC.
We show that these changes in cFTLE ridges can be used to understand the change in value functions useful for optimal control.
We also show that cFTLE ridges can be a powerful tool for numerically identifying switching surfaces and approximations to time-varying invariant manifolds.
This could potentially aid in using various techniques from non-smooth analysis, and knowledge about global bifurcations to analyze and predict the oncoming changes in system stability as hyperparameters are varied.

\subsection{Interpretations of cFTLE}
We first discuss the interpretations of cFTLE ridges in comparison with the traditional view of FTLE when applied to the transport of agents through a flow field. 
We will vary the relative cost of actuation versus state error, given by the ratio $R/Q$, the time horizon $T_H$, position of the goal and analyze how this impacts the cFTLE. In the limit of an infinite cost of control, the cFTLE should converge to the traditional passive FTLE. 

\subsubsection{Transport Barriers}
Forward time cFTLE ridges form transport barriers similar to traditional FTLE ridges~\cite{wilson2009lagrangian}.
This can be observed from Figure~\ref{fig:regimes}, where we have two patches of agents moving towards a common goal.
We observe that in the time of advection, the purple patch and yellow patch are separated by the cFTLE ridge as a dividing barrier.
Agents that fall on a cFTLE ridge are further stretched apart at later times, which will be especially highlighted later in Section~\ref{sec:RL}.
Such analysis is useful in the context of analyzing control laws for large-scale transport of under-actuated mobile sensors in fluid flows. For instance, determining the deployment locations for intelligent mobile sensors is a challenging task with significant cost and schedule complications~\cite{Talley:19a}. cFTLE can be potentially useful in finding locations where the control law is able to transport agents effectively.

\subsubsection{Rendezvous Through Lobe Transport}
\begin{figure}[t]
    \centerline{\includegraphics[scale=.3,trim = {0 60mm 0 40mm },clip,]{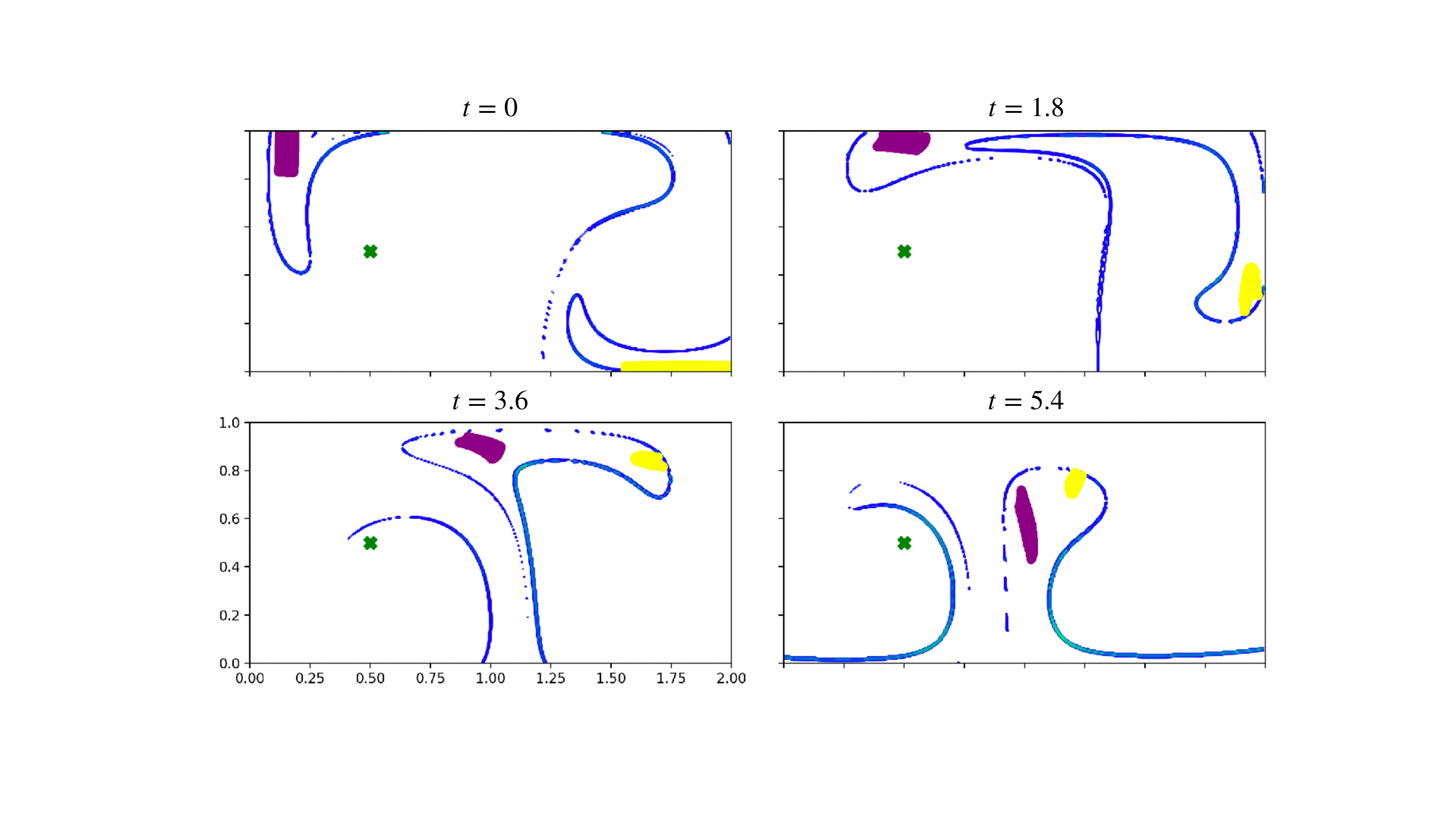}}
    \caption{cFTLE boundaries enclose regions in the domain which accrue similar cost. Despite the fact that the purple and the yellow patch start far from each other, they meet in short time since they are in the same lobe. This can be viewed particularly at $t = 3.6$ on, where the yellow and purple patch are within the same bounded cFTLE region. These structures are not apparent when viewing the control law or the flow field alone. An application of visualizing such regions in the domain may aid in solving rendezvous problems in the ocean where two distinct groups of agents may need to meet in finite time}
    \label{fig:rendezvouz}
\end{figure}

Building upon the idea of transport barriers from the previous section, cFTLE barriers can surround agents with barriers that persistently keep the agents trapped.
This behavior has been exhaustively studied in dynamical systems theory~\cite{wiggins2013chaotic, rom1990analytical}.
Attracting and repelling invariant manifolds intersect and form lobes that trap passive particles. 
Traditionally, attracting and repelling FTLE ridges have been used to visualize such trapping regions.
We observe this in the context of cFTLE in Figure~\ref{fig:rendezvouz}. 
Despite not plotting the attracting cFTLE ridges, we can see that the repelling cFTLE barriers, inside of which the purple and yellow patches are placed, form a single trapping structure that brings the two patches to meet in forward time, while they move towards the goal.
This idea is potentially useful in the context of multi-agent planning, where drop locations can be found that allow agents to rendezvous in forward time. Such rendezvous problems are important in the context of path planning~\cite{song2017multi}.

\subsection{Connections to Optimal Control}

\begin{figure}[t]
    \centerline{\includegraphics[scale = .27,trim = {0 120mm 0 120mm },clip,]{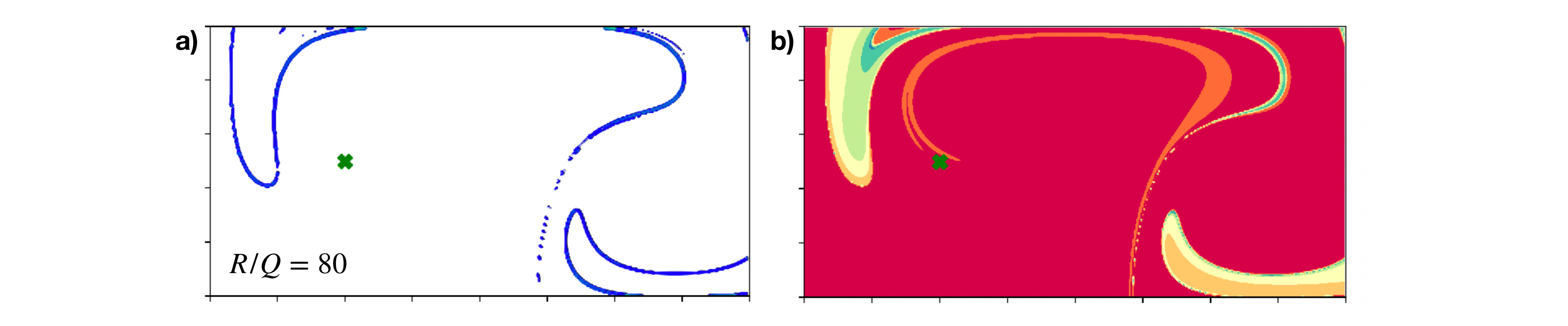}}
    \caption{The plot (a) shows the control FTLE ridges for double gyre setup at an $R/Q$ ratio of 80 on the left, with a time of advection $T_A = 15, T_H = 3$. The agents are tracking the goal set point at $\bx_{goal} = [0.5, 0.5]$. On the right (b), we show the accumulated state error $J_F$ for a meshgrid of agents in the domain for the same $R/Q$ and $T_H$. We see that the cFTLE ridges are located where there is a drastic change in color on the right figure. This indicates that the cFTLE ridges at a particular $T_A$ partition the domain by the cost spent (i.e., large change in $J_F$).}
    \label{fig:j_vs_ridge}
\end{figure}

\begin{figure}[t]
    \centerline{\includegraphics[scale = .27,trim = {0 120mm 0 90mm },clip,]{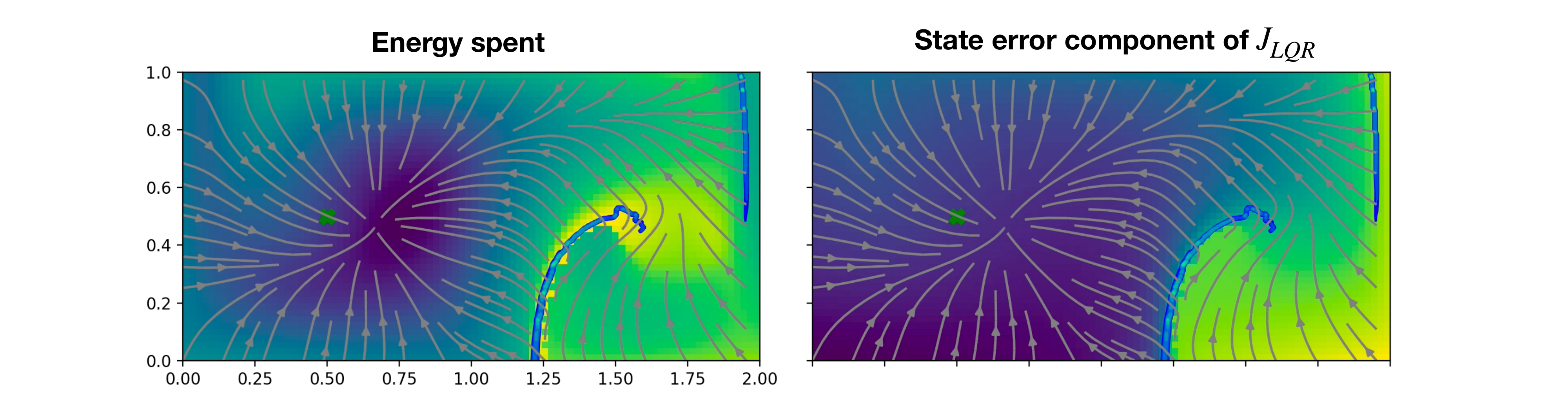}}
    \caption{On the left, we plot the instantaneous energy spent by an agent starting from a meshgrid of initial conditions moving to a goal location on the left ($\bx_{goal} = [0.5, 0.5]$) at $t = 0$. On the right, we plot the accumulated state error over the time horizon of MPC ($T_H = 4.5, R/Q = 15$). In both plots, the streams are the control law generated by all the MPC agents, and the cFTLE ridges are shown where the time of advection is equal to the time horizon of the MPC, $T_A = T_H$. We observe here that the cFTLE ridge overlaps with regions where actuation energy expenditure peaks locally. Also, similar to Figure \ref{fig:j_vs_ridge}, the cFTLE divides the domain where the cost gradient is the largest.}
    \label{fig:j_vs_ridge_MPC}
\end{figure}

In this section, we discuss cFTLE from the perspective of optimal control. Inspired by the first term in Eq.~\eqref{eq:cost_fn1}, Figure~\ref{fig:j_vs_ridge} (right) displays the value function landscape of 
\begin{equation}
    J_F = \|\bx(T_A) - \bx_{goal} \|_2^2,
    \label{eq:J_F}
\end{equation}
 which, intuitively, is the square of the agent distance from the goal after a time $T_A$ has elapsed, as a function of where the agent initially started.
By definition, $J_F$ only depends on the final state error at time $T_A$, in contrast to $J_{LQR}$, which depends on the state error and control energy spent integrated over the trajectory.
The general idea of this section is that we will use $J_F$ as a simple model to gain insight into $J_{LQR}$.

When $J_F$ is viewed in comparison to the cFTLE on the left in Figure~\ref{fig:j_vs_ridge}, we observe that the cFTLE ridges overlay with regions where there is a sharpest change in value of $J_F$. 
This indicates that cFTLE ridges coincide with regions of high sensitivity in the state error of a cost function, which will be further verified later in this section.
This is in agreement with the definition of the FTLE, which is a measure of sensitivity of final conditions to infinitisimilly close initial conditions.

Furthermore, the left plot in Figure~\ref{fig:j_vs_ridge_MPC} displays the energy spent by the control action taken, $\|\bu^*\|$, as a function of initial condition.
We see that cFTLE ridges fall on regions where the energy spent by the agent is the largest. 
We know that Eq.~\eqref{eq:ctrl_grad} implies that $\|\bu^*\|_2^2 = \|\ \mathbf{R}^{-1} \nabla J_{LQR}\|^2_2$, and, we established in the previous paragraph that cFTLE ridges are connected to sharpest \emph{change} in $J_F$. 
In image processing, the sharpest change of an image or function value (e.g., $J_F$) can be found by taking the magnitude of the gradient of the function ($\|\nabla J_F\|_2^2$). 
These results specifically highlight a connection between cFTLE ridges, sensitivity in value functions $\|\nabla J\|_2^2$, and the instantaneous energy spent by optimal control policies $\| \bu^* \|_2^2$.

It is also possible to further explore this connection mathematically. When using the fact that $\bx({T_A})= \boldsymbol{\hat\Phi}_{t_0}^{t_0+T_A}(\bx(t_0))$, the gradient of $J_F$ is given by 
\begin{equation}
    \nabla J_F = 2[\mathbf{x}(T_A) -\mathbf{x}_{goal}] ^T \mathbf{D}\boldsymbol{\hat\Phi}_{t_0}^{t_0 +T_A}.
    \label{eq:grad_jf}    
\end{equation}
By comparing to the FTLE definition, we can observe that the flow map operator used in the computation of FTLE is deeply connected to the gradient of cost functions with quadratic state error. Specifically, the maximum cost function gradient is related to the maximum singular value of the flow map Jacobian, which is related to the FTLE.  

At each instant over the horizon, we compute a complete finite-horizon optimization and take only the first control action, a standard practice of MPC. 
On the right plot in Figure~\ref{fig:j_vs_ridge_MPC}, we plot the integrated state error cost of this finite-horizon optimization at the first step - the first term of $J_{LQR}$.
If we were to make $T_A = T_H$, the $J_F$ plot in Figure~\ref{fig:j_vs_ridge} would look similar to the right plot in Figure~\ref{fig:j_vs_ridge_MPC}, further highlighting the connection between cFTLE and sensitivity in $J_{LQR}$.
We observed across our simulations that the cost function from $J_F$ and $J_{LQR}$ had very similar qualitative features, and one could potentially use whichever is easier to compute to gain insight into the other.

In model-free control approaches, such as reinforcement learning, one is interested in generating fast \emph{approximations} to the value function as the agent is moving through the flow field. 
These results show that cFTLE ridges can have a potential impact on estimating these value functions when such model free approaches are used. We will further explore reinforcement learning in Section~\ref{sec:RL}.

\subsection{Deformation of cFTLE with change in MPC parameters}
We have now shown that cFTLE ridges contain information about the sensitivity of a control policy, the regions of large magnitude of control, and the boundaries between finite-time invariant sets. 
In this section, we study how these ridges depend on the aggressiveness of the control law and the location of the goal. 
Particularly, these two parameters direct the strength and magnitude of fluxes generated by the model predictive control policy over the unsteady flow field.
We also study change in cFTLE as the time horizon $T_H$ is varied. 
This parameter changes the spatial complexity of the control policy.
These results can be potentially useful for anticipating changes in the cFTLE for a different parameter value. 
In the next section (Section~\ref{sec:RL}), we finally use the ideas learned in this section to interpret a policy generated through another popular optimization-based control approach -- reinforcement learning.

\subsubsection{Varying the Cost of Control, $R/Q$}
\begin{figure}[t]
    \centerline{\includegraphics[scale=.28,trim = {0 55mm 0 63mm },clip,]{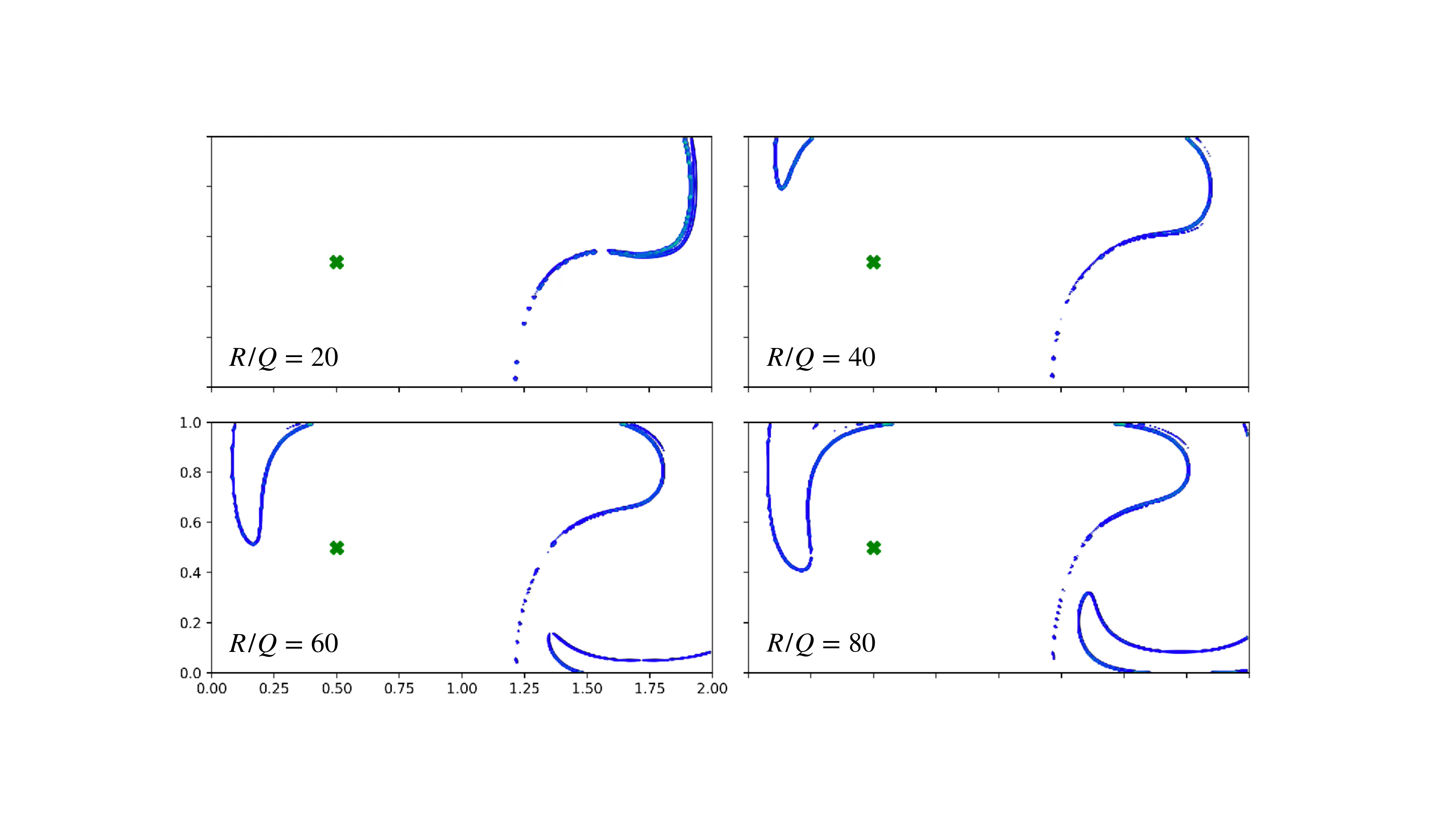}}
    \caption{This figure shows the cFTLE ridges for four different $R/Q$ ratios of ${20, 40, 60, 80} $, and how the barriers shrink as the aggressiveness of the MPC increases. Here, the time of advection for cFTLE $T_A = 15$, and the time horizon of MPC is $T_H = 3$. We observe that as the aggressiveness of the control decreases, the cFTLE approaches the passive FTLE.}
    \label{fig:cFTLE_vs_rq}
\end{figure}

Figure~\ref{fig:cFTLE_vs_rq} shows cFTLE ridges plotted for incremental changes in $R/Q$ ratio. 
We find that as the $R/Q$ ratio increases, the aggressiveness of the MPC strategy decreases, and the average magnitude spent for control decreases, while the direction of control for each state does not change appreciably. 
Extremely aggressive control (small $R/Q$) pushes the cFTLE ridges further towards the right when the goal is in the left gyre, and shrinks the size of the left cFTLE ridge completely. 
For less aggressive control (large $R/Q$), the cFTLE ridges lengthen and begin to approach passive FTLE ridges. 
These results are intuitive, given that when the control becomes less aggressive, the agent behaves more like a passive particle due to lack of control authority.

\subsubsection{Changing the Goal Location, $\bx_{goal}$}

\begin{figure}[t]
    \centerline{\includegraphics[scale=.26,trim = {0 110mm 0 110mm },clip,]{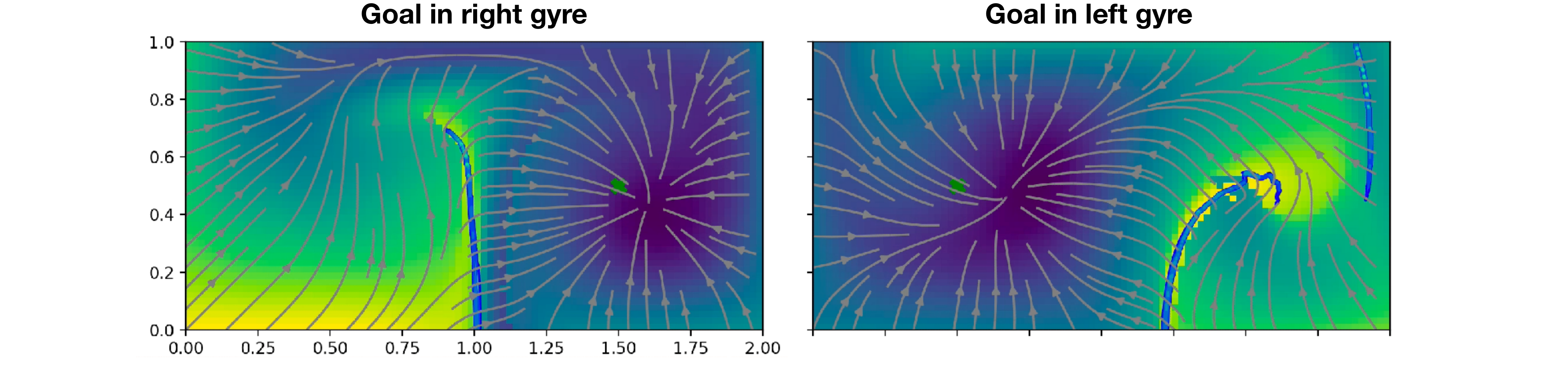}}
    \caption{This figure highlights the change in cFTLE ridges as a function of the goal location. Here, $R/Q = 15, T_A = T_H = 4.5$. The streamplot on both plots show the control law. The color in the background shows the regions where most energy is spent (brighter yellow shows larger energy spent and darker shows less). On the left plot, the goal location is set to the right gyre, and on the right plot, the goal location is set to the left gyre. We observe how the fluxes generated by the control law moves the passive FTLE ridge. A left flux moves the cFTLE ridge to the right, and vice versa. This figure continues to highlight the connection between the terms in the cost function and cFTLE ridge as in Figure~\ref{fig:j_vs_ridge}.}
    \label{fig:cFTLE_vs_goal}
\end{figure}

The results of changing the goal location can be found in Figure~\ref{fig:cFTLE_vs_goal}. 
When a control law or policy acts over an unsteady flow field, it generates an added flux to transport particles towards the goal. 
This is visualized by the grey arrows in Figure~\ref{fig:cFTLE_vs_goal}. 
Here, we observe how the cFTLE changes depending on this newly added flux. 
For instance, creating a sink in the right gyre moves the cFTLE ridge to the left from the initial position of the corresponding passive FTLE ridge. 
Creating a sink to the left moves the cFTLE ridge to right. 
This shows us that cFTLE ridges move in the direction opposite to the flux of flow generated.

\subsubsection{Varying the Time Horizon, $T_H$}

\begin{figure}[t]
    \centerline{\includegraphics[scale=.26,trim = {0 20mm 0 40mm },clip,]{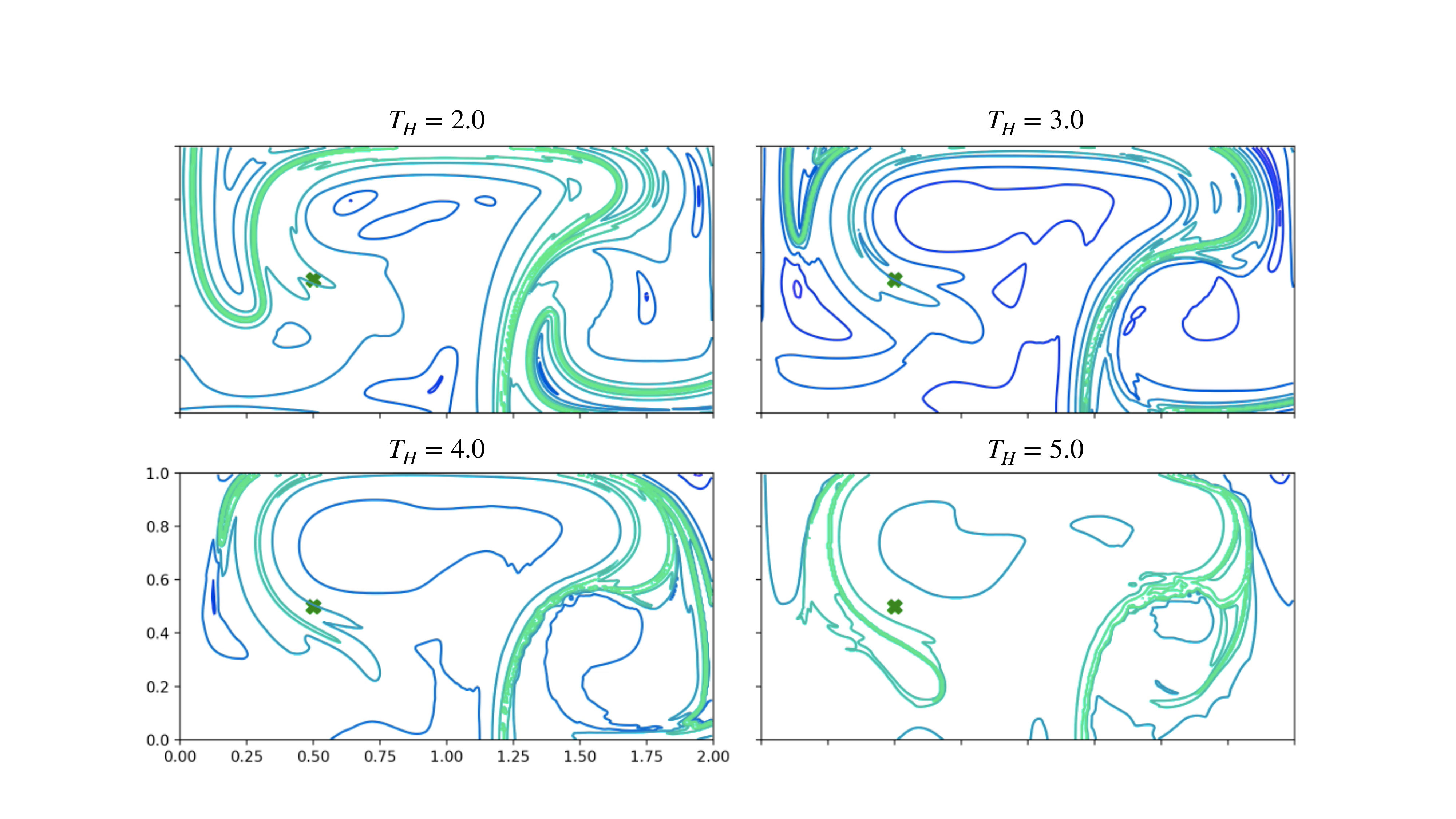}}
    \caption{This figure shows the change in cFTLE field for time horizon values $T_H = \{2, 3, 4, 5\}$, using an advection time of $T_A = 15$ and $R/Q = 50$ at $t_0 = 0$. Unlike in Figure~\ref{fig:cFTLE_vs_rq}, the cFTLE ridge changes due to more intelligent use of control and not because of sheer use of greater effort. The cFTLE ridges exhibits more structures, and the curves have more branches due to several ridges coming close to each other and collapsing on each other.}
    \label{fig:th_vary}
\end{figure}

\begin{figure}[t]
    \centerline{\includegraphics[scale =.26,trim = {0 95mm 0 130mm },clip,]{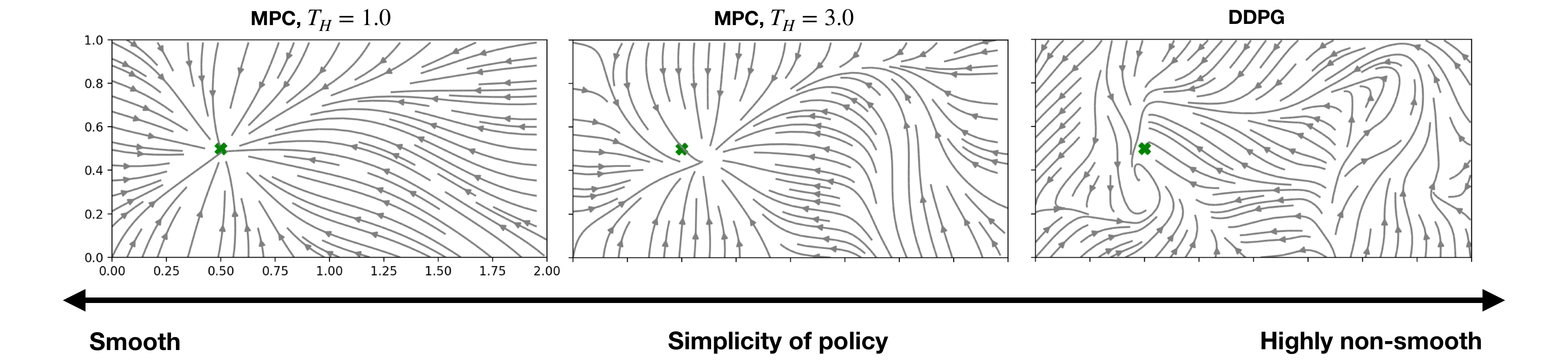}}
    \caption{This figure shows the range of policies that can be generated at time $t = 0$. The policies considered in this paper vary with time $t$. The left and middle policies are generated using MPC with a time horizon of $1.0$ and $3.0$, respectively. The policy on the right is generated by using DDPG and can be considered the ``infinite-horizon'' case. All policies have an $R/Q$ ratio of $70$. From left to right, there is an increase in spatial complexity of the policy as it begins to use more information of the flow field. This spatial complexity in turn influences the cFTLE ridges computed, which can be seen in Figure~\ref{fig:th_vary} and Figure~\ref{fig:RL_cFTLE}.}
    \label{fig:policies}
\end{figure}

Figure~\ref{fig:th_vary} shows the change in cFTLE field as the time horizon parameter $T_H$ is varied. 
Note that we plot the contours of the FTLE field as opposed to the FTLE ridges since agents following the policy cause the ridges to be less sharp, and therefore, more challenging to extract and visualize through thresholding.
Unlike the deformation of the cFTLE in the previous section, looking at the deformation of cFTLE in Figure~\ref{fig:th_vary}, we observe that cFTLE undergoes a deflation-type effect and the cFTLE field becomes non-smooth with several branches.
In Figure~\ref{fig:policies}, the first two figures show the policy changing as a function of $T_H$. 
The time horizon parameter captures the amount of future knowledge of the flow field being incorporated in taking the present action. 
Therefore, the policy changes from a naive sink-like behavior to become spatially more complex as the time horizon is increased. 
These spatially complex policies make better use of actuation to escape regions or lobes dictated by the \emph{passive} FTLE, where the agent will be advected away from the goal and increase the cost function value.

We also observe in Figures~\ref{fig:j_vs_ridge_MPC}~and~\ref{fig:cFTLE_vs_goal} that the cFTLE highlights regions in space where the policy becomes spatially discontinuous. 
The streamlines on either side of the cFTLE ridge are in different directions.
This is intuitive, as regions spatial discontinuity in a vector field can potentially also generate large Lyapunov exponent values, in addition to effects from shear and normal hyperbolicity.
When using optimal control methods, the policy can often make the controlled system a non-smooth dynamical system~\cite{sastry2013nonlinear}.
These results shows that the cFTLE can be used for switching manifold detection in controlled systems.

\section{Example from Reinforcement Learning}\label{sec:RL}

\begin{figure}[t]
    \centerline{\includegraphics[scale = .27,trim = {0 50mm 0 30mm },clip,]{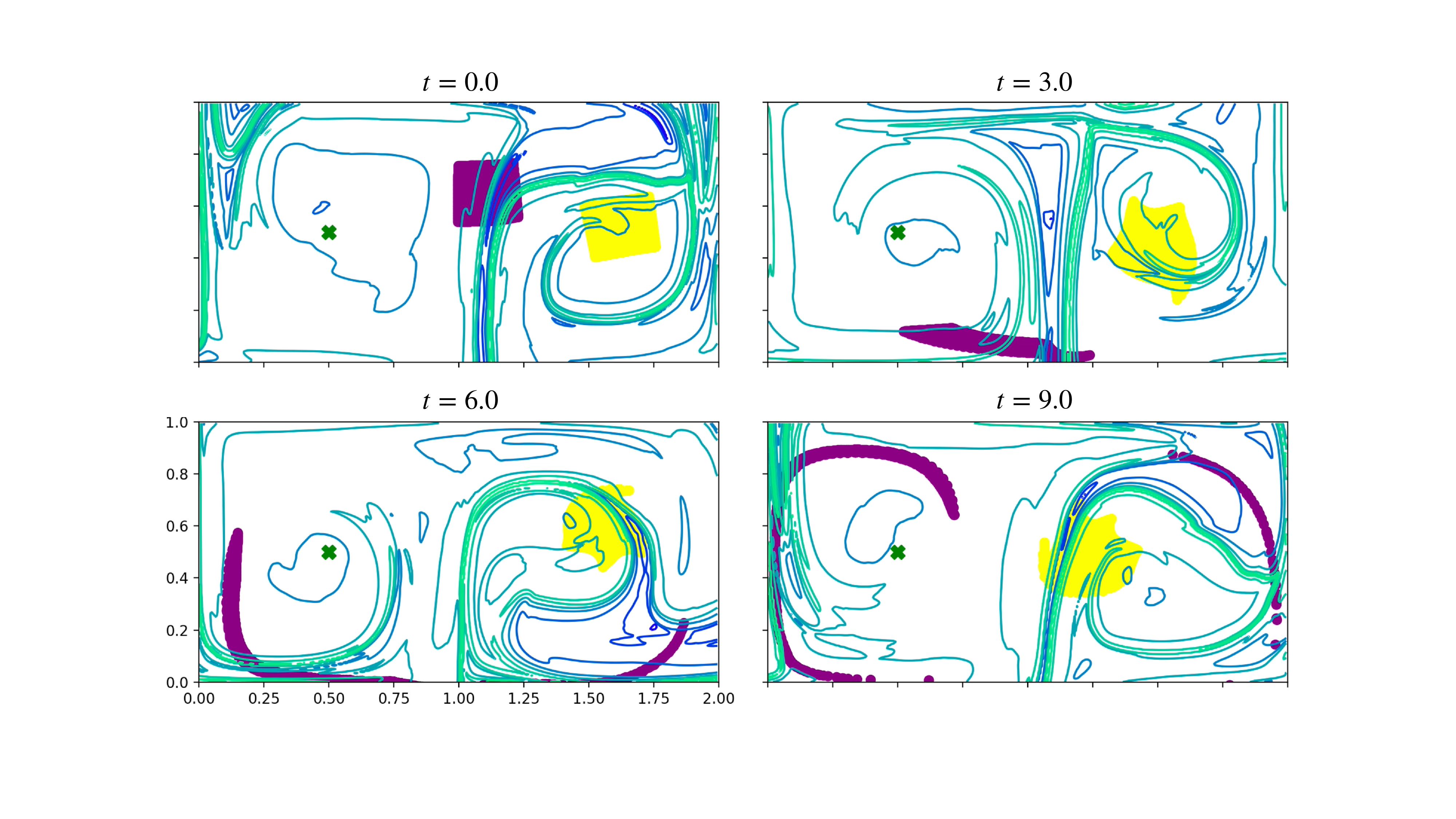}}
    \caption{This example shows the repelling cFTLE field plotted for a policy generated using DDPG with $R/Q = 70$. Much like Figure~\ref{fig:regimes}, we use two patches to highlight transport mechanisms. We also use the same initial patch locations for comparison.}
    \label{fig:RL_cFTLE}
\end{figure}

Thus far, we have discussed the use of cFTLE exclusively on policies generated from model predictive control.  
However, cFTLE can be used to analyze policies generated from any control or planning methods.
Reinforcement learning (RL) is a powerful alternative paradigm to model predictive control for policy generation~\cite{recht2019tour,Brunton2022book}. 
A key difference between the two paradigms is that RL is typically model free, i.e., policies can be generated without any knowledge of the flow field, but rather from data alone.
In this study, we use the deep deterministic policy gradient (DDPG) method~\cite{lillicrap2015continuous} available in the stable baselines package~\cite{stable-baselines3} to demonstrate the use of cFTLE analysis on reinforcement learning policies.

Much like model predictive control, RL requires the definition of a reward function. 
Typically, the cumulative reward over the simulation is computed.
For this, we make use of the negative of the discretized LQR cost function, defined in Eq.~\eqref{eq:cost_fn1} due to the fact that the convention is to \emph{maximize} the reward, as opposed to minimizing the cost in methods such as MPC.
We make use of an Euler time-stepping scheme to propagate an agent's state forward in time with a time step of $0.1$. 
Each complete RL simulation (or episode) is run for $800$ time steps, which totals to 80 units as in the MPC case.

At a high level view, RL runs multiple episodes with random policies and random initial states to generate data of cumulative rewards over each episode.
This data is then used to iterate towards a policy that maximizes the reward over each episode, also known as the cumulative reward. 
The setup of DDPG used in this paper does not explicitly incorporate a time horizon; however, the resulting policy can be interpreted as a solution of an infinite horizon problem.
The right plot in Figure~\ref{fig:policies} shows an RL policy generated over the double gyre flow field for the same objective outlined in the problem setup section.
Metrics such as the average cumulative reward are often used to judge the performance of the policy~\cite{sutton2018reinforcement}.
These metrics do not provide a spatial description or boundaries in space where the reward declines, which is possible through cFTLE analysis.
RL generates a lookup table policy function, described by a neural network $\bu(\bx, t)$ for each state in space and time.
Therefore, unlike in the MPC case, computing the FTLE of an RL policy does not require the use of interpolation as mentioned earlier.
As in the previous section, we plot the cFTLE \emph{field} as opposed to cFTLE ridges to better visualize these invariant structures.

Figure~\ref{fig:RL_cFTLE} shows the integration of two patches of different colors with the same initial configuration of Figure~\ref{fig:regimes}. 
We observe in contrast to Figure~\ref{fig:regimes}, that, since the purple patch starts on the cFTLE ridge as opposed to on one side of it, it undergoes large deformation normally along the cFTLE ridge in the middle of the domain, while the yellow patch maintains cohesion for longer time. 
However, at time $t=9.0$, we see that a part of the yellow patch lies on the the other side of the cFTLE ridge, therefore, we can predict that in the future, the yellow patch will also be stretched out.
Since the RL policy can be viewed as an infinite-horizon extension of the MPC case, we can interpret the several streaks of large cFTLE value as passive FTLE structures that have collapsed very close to each other. 
This can be visualized at $t=3.0$ where the FTLE ridge takes a ``P'' shape that appears as a branch off the middle hyperbolic LCS, which occurs due to the shrinking and collapsing of a passive FTLE lobe as seen from the large time horizon case.
This is similar to the case of time horizon $4.0$ and $5.0$ in Figure~\ref{fig:th_vary}.

Value functions for RL methods can find interesting pathways through passive dynamics that traditional methods miss. 
This is due to the highly non-convex nature of RL algorithms, which makes it challenging to bound and analyze the policies generated through these methods.
cFTLE provides a principled method of analyzing these policies. 
Moreover, given that value functions are connected to the cFTLE as demonstrated earlier, we can potentially use cFTLE to minimize the steps required to compute policies, leading to faster convergence.

\section{Discussion \& Conclusion} 
In this paper, we studied the use of the FTLE method on active agents that use actuation to move towards a fixed goal in an unsteady flow field. 
In particular, we focused on the FTLE analysis of agents using model predictive control and reinforcement learning independently to generate policies in a double gyre flow field. 
Broadly speaking, this work can be interpreted through the lens of its potential applications to the following problems:

\emph{1) Invariant sets}: We first compute the cFTLE for agents using policies generated through MPC. 
We confirm that, much like passive FTLE, cFTLE can be used to find barriers that separate different regions of agents following a particular policy. 
This suggests that cFTLE can be useful to identify finite-time invariant sets and invariant manifolds for active agents, which is particularly useful in the context of navigating unsteady aperiodic flow fields as coherent structures only persist for finite time. 
Invariant sets are useful in optimal control from the perspective of analyzing robustness and designing model predictive control policies~\cite{rakovic2007optimized}.
The areas with smaller cFTLE values correspond to regions where the control is robust to small changes in the parameters of the controller.

\emph{2) Multi-agent path planning}: Next, we find that cFTLE can be used to identify trapping regions or regions closed off by barriers or lobes.
When viewing these trapping regions at large scales in an unsteady flow field, we can find regions in which, when following the optimal policy, agents can move towards a goal cohesively without getting separated at large distances.
When viewing these barriers at small scales, they can be used to find deployment locations where the agents would not collide with physical objects in the flow field or other agents following the same policy.

\emph{3) Optimal control theory}: Value functions are scalar functions defined over the domain that map each state to the expected sum of future rewards from that state. 
From this perspective, control policies (in the context of the kinematic models considered in this paper) are gradients of the value function that direct individual agents to ``climb up'' to states with greater values.
Value functions are, generally speaking, solutions to the HJB equation.
Modern equation-free control methods such as reinforcement learning ultimately attempt to solve the HJB equation approximately.
Recent works along this direction in the context of navigating flow fields include~\cite{wiggert2022navigating,doshi2023energy}.
In our work, we show that cFTLE ridges can be deeply connected to the boundaries of a value function, where the cost function is most sensitive to perturbations. 
We also highlight that regions most sensitive to perturbation are also the regions where greatest energy is spent by the agent.
These results can potentially aid developing methods for faster computation of value functions or policies. 
For example, only the value of states where the cFTLE field is large can be updated, since other states will be less sensitive to perturbations when a parameter like the $R/Q$ is varied.
This can save on redundant computation at states which do not change value under perturbation.
Furthermore, a deeper connection can be potentially highlighted mathematically connecting the singular values/vectors of the flow map Jacobian, which are related to cFTLE, and the value function. 

\emph{4) Bifurcation and stability analysis}: When the $R/Q$ ratio approaches 0, the policy has the actuation capability to turn the goal state into an attracting fixed point.
As we increase the $R/Q$ ratio, the policy loses the ability to create a strongly attracting fixed point, and instead, creates a limit cycle oscillation around the goal state, as was reported previously in~\cite{krishna2022finite}.
When there is no actuation (i.e., in the limit of $R/Q$ going to $\infty$), trajectories in the double gyre are chaotic.
This points to the possibility that changing the $R/Q$ ratio generates successive bifurcations, ending in chaos.
Invariant manifolds, cFTLE ridges in this case, play a crucial role in global bifurcations, which could be useful in understanding the change in stability described above~\cite{guckenheimer2004fast,guckenheimer2015invariant}.
In the context of transport of autonomous agents, this would allow us to understand and predict the onset of loss of stability around the goal for different $R/Q$ ratios, extrapolating from known positions of the cFTLE ridges.
This would ultimately aid in tuning the $R/Q$ parameter.

The latter part of the paper focuses on the deformation of cFTLE ridges, or rather, the intuition behind how the cFTLE ridges move when new fluxes from the control policy are added to the passive unsteady flow field.
We find that the cFTLE ridge moves \emph{opposite} to the direction of the added flux from the flow field.
We also find that control can cause a trapping barriers to shrink and collapse into multiple branches emerging from a larger ridge.
This can highlight regions (or lobes) in space that are integral in escaping or entering to move towards the goal effectively.
An interesting observation is that cFTLE ridges can highlight the presence of switching manifolds in controlled systems, where the optimal policy renders the controlled dynamical system non-smooth.
This is possible since regions of discontinuity in  a vector field can generate large Lyapunov exponents.
Finally, the computation of cFTLE ridges itself is challenging. Control policies generate sources and sinks, where innovations in the computation of compressible FTLE can be used~\cite{gonzalez2016finite}.
Potentially, many of these ideas can be tested in realistic experimental systems with three-dimensional flows and multi-scale turbulence.

\section*{Acknowledgments}
The authors acknowledge funding support from the National Science Foundation AI Institute in Dynamic Systems (grant number 2112085) and the US Air Force Office of Scientific Research (FA9550-21-1-0178). 

\appendix

\section{Attracting cFTLE}
\begin{figure}[t]
    \centerline{\includegraphics[scale=.27,trim = {0 40mm 0 30mm },clip,]{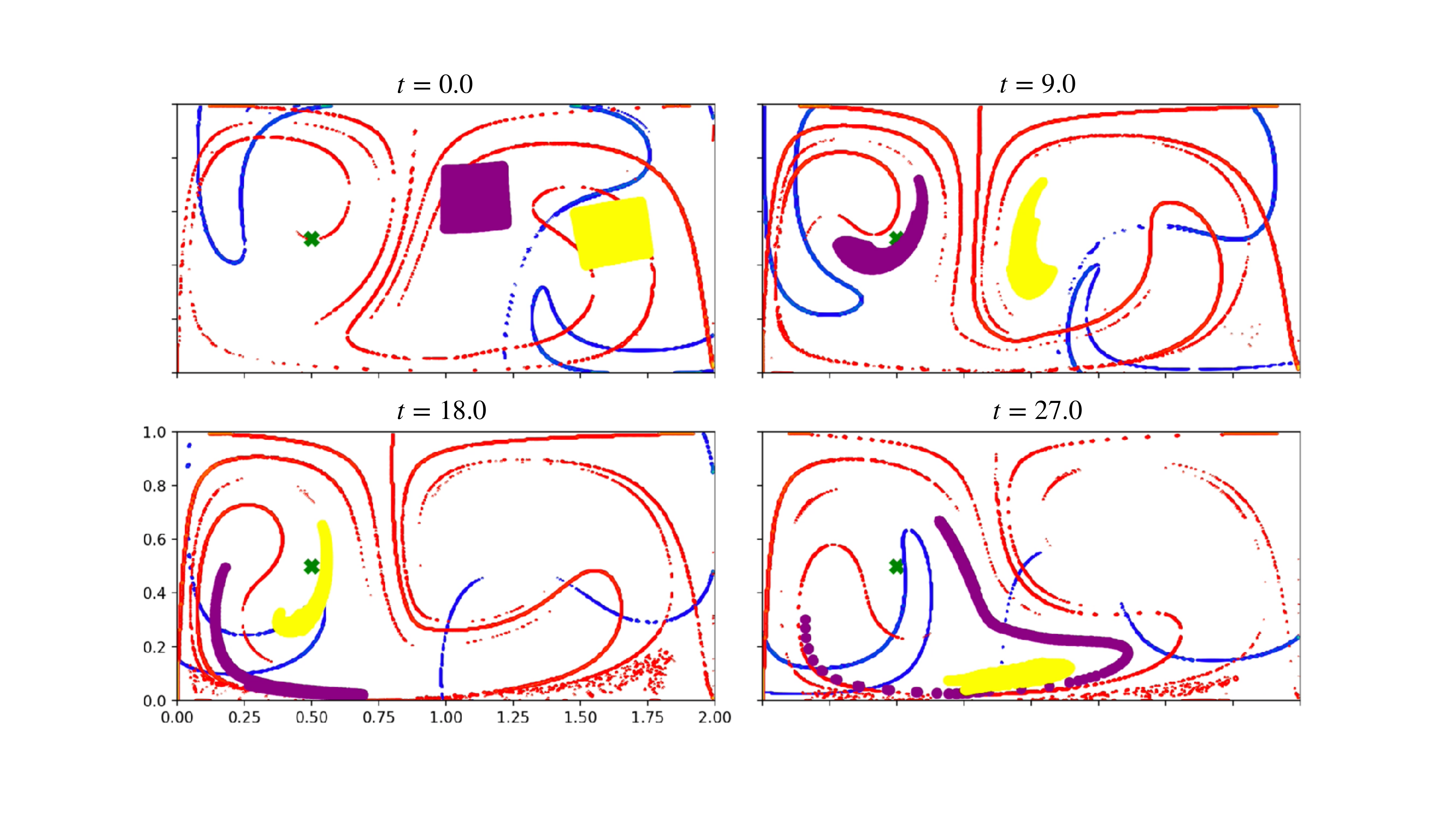}}
    \caption{This figure shows two patches of agents - one purple, and one yellow starting on opposite sides of a repelling cFTLE ridge computed with $T_A = 15$ and evolving through the unsteady flow field with control. We also plot the attracting cFTLE ridges to show the full picture of high the cFTLE ridges govern the dynamics. We see that the attracting cFTLE ridges highlight curves of long term attraction}
    \label{fig:attract}
\end{figure}
In this paper, much of the discussion has revolved around the repelling cFTLE. 
This is due to fact that the policy for moving towards a goal is well known in forward time. 
To compute an attracting cFTLE, it is important to know the policy in backward time.
This is easier if the flow field is periodic, since the policy can simply be reversed in time sequence.
However, in more realistic situations, it is challenging to have the flow field data and a control policy at times \emph{before} the initial condition, which is generally not possible.
Therefore, the attracting cFTLE can only be visualized at times $(t - t_0) > T_A$, where $t_0$ is the time point at which the flow field data is available.

\begin{figure}[t]
    \centerline{\includegraphics[scale=.27,trim = {0 120mm 0 110mm },clip,]{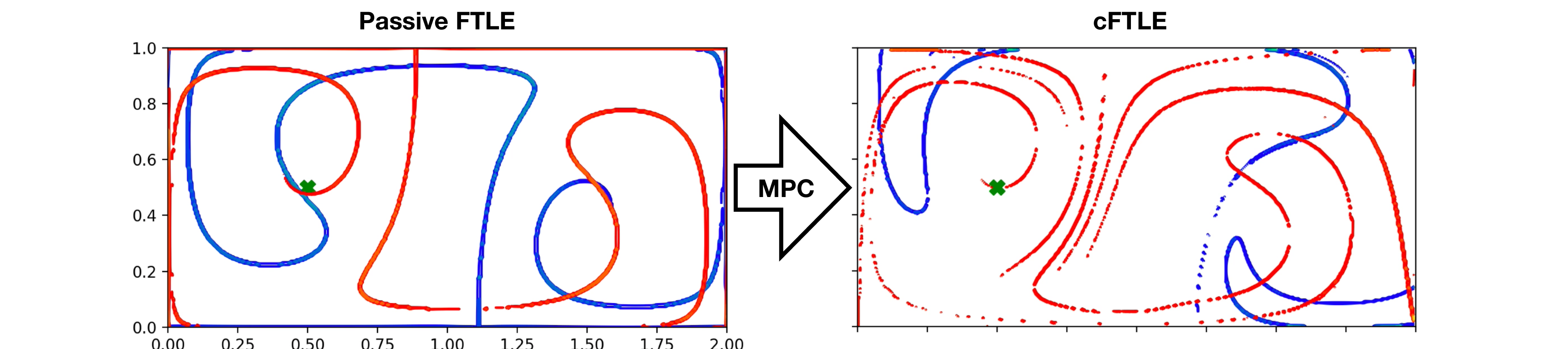}}
    \caption{We plot the passive FTLE ridges in the uncontrolled double gyre on the left, and the deformed ridges under the action of MPC at $T_H = 3.0, R/Q = 80$ on the right. The blue curves are repelling structures, and the red curves are attracting structures.}
    \label{fig:pass_ctrl}
\end{figure}

In Figure~\ref{fig:attract}, the data from Figure~\ref{fig:regimes} is plotted, but this time with the attracting cFTLE also visualized. 
The attracting cFTLE looks markedly different from the passive FTLE in Figure~\ref{fig:pass_ctrl}. 
The double gyre system contains six fixed points counting the corners of the domain and the two on the edges. 
When control is applied, the unstable manifolds from each of these fixed points connect to an attractor formed around the goal by the control policy. 
This causes the attracting cFTLE to have multiple curves leading into the domain, taking a spiral structure.
In the passive case, these unstable manifolds align with the boundaries of the domain. 
Depending on the aggressiveness of the control policy due to the hyperparameters, the attracting set can either be a fixed point, a limit cycle, or a chaotic set. 
In Figure~\ref{fig:attract}, $T_A$ is much shorter in timescale than the time it takes for all the trajectories to fall on the attractor. 
This can lead to the attractor not being fully resolved when the cFTLE ridges are plotted.
Being able to visualize attracting sets is useful in the context of global bifurcation analysis~\cite{guckenheimer2015invariant}.
We observe that, the attracting cFTLE forms curves in space, to which, patches latch onto in forward time.

 \small{
 \setlength{\bibsep}{4.5pt}
 \bibliographystyle{IEEEtran}
 \bibliography{references}

\begin{thebibliography}{10}
\providecommand{\url}[1]{#1}
\csname url@samestyle\endcsname
\providecommand{\newblock}{\relax}
\providecommand{\bibinfo}[2]{#2}
\providecommand{\BIBentrySTDinterwordspacing}{\spaceskip=0pt\relax}
\providecommand{\BIBentryALTinterwordstretchfactor}{4}
\providecommand{\BIBentryALTinterwordspacing}{\spaceskip=\fontdimen2\font plus
\BIBentryALTinterwordstretchfactor\fontdimen3\font minus
  \fontdimen4\font\relax}
\providecommand{\BIBforeignlanguage}[2]{{%
\expandafter\ifx\csname l@#1\endcsname\relax
\typeout{** WARNING: IEEEtran.bst: No hyphenation pattern has been}%
\typeout{** loaded for the language `#1'. Using the pattern for}%
\typeout{** the default language instead.}%
\else
\language=\csname l@#1\endcsname
\fi
#2}}
\providecommand{\BIBdecl}{\relax}
\BIBdecl

\bibitem{fossum2019toward}
T.~O. Fossum, G.~M. Fragoso, E.~J. Davies, J.~E. Ullgren, R.~Mendes,
  G.~Johnsen, I.~Ellingsen, J.~Eidsvik, M.~Ludvigsen, and K.~Rajan, ``Toward
  adaptive robotic sampling of phytoplankton in the coastal ocean,''
  \emph{Science Robotics}, vol.~4, no.~27, 2019.

\bibitem{chai2020monitoring}
F.~Chai, K.~S. Johnson, H.~Claustre, X.~Xing, Y.~Wang, E.~Boss, S.~Riser,
  K.~Fennel, O.~Schofield, and A.~Sutton, ``Monitoring ocean biogeochemistry
  with autonomous platforms,'' \emph{Nature Reviews Earth \& Environment},
  vol.~1, no.~6, pp. 315--326, 2020.

\bibitem{zhang2021system}
Y.~Zhang, J.~P. Ryan, B.~W. Hobson, B.~Kieft, A.~Romano, B.~Barone, C.~M.
  Preston, B.~Roman, B.-Y. Raanan, D.~Pargett \emph{et~al.}, ``A system of
  coordinated autonomous robots for {L}agrangian studies of microbes in the
  oceanic deep chlorophyll maximum,'' \emph{Science Robotics}, vol.~6, no.~50,
  2021.

\bibitem{bellingham2007robotics}
J.~G. Bellingham and K.~Rajan, ``Robotics in remote and hostile environments,''
  \emph{Science}, vol. 318, no. 5853, pp. 1098--1102, 2007.

\bibitem{wynn2014autonomous}
R.~B. Wynn, V.~A. Huvenne, T.~P. {Le Bas}, B.~J. Murton, D.~P. Connelly, B.~J.
  Bett, H.~A. Ruhl, K.~J. Morris, J.~Peakall, D.~R. Parsons, E.~J. Sumner,
  S.~E. Darby, R.~M. Dorrell, and J.~E. Hunt, ``Autonomous underwater vehicles
  ({AUVs}): Their past, present and future contributions to the advancement of
  marine geoscience,'' \emph{Marine Geology}, vol. 352, pp. 451--468, 2014.

\bibitem{rhoads2013minimum}
B.~Rhoads, I.~Mezi{\'c}, and A.~C. Poje, ``Minimum time heading control of
  underpowered vehicles in time-varying ocean currents,'' \emph{Ocean
  Engineering}, vol.~66, pp. 12--31, 2013.

\bibitem{song2017multi}
Z.~Song, D.~Lipinski, and K.~Mohseni, ``Multi-vehicle cooperation and nearly
  fuel-optimal flock guidance in strong background flows,'' \emph{Ocean
  Engineering}, vol. 141, pp. 388--404, 2017.

\bibitem{song2015anisotropic}
Z.~Song and K.~Mohseni, ``Anisotropic active {L}agrangian particle swarm
  control in a meandering jet,'' in \emph{2015 54th IEEE Conference on Decision
  and Control (CDC)}.\hskip 1em plus 0.5em minus 0.4em\relax IEEE, 2015, pp.
  240--245.

\bibitem{shadden2005physd}
S.~C. Shadden, F.~Lekien, and J.~E. Marsden, ``Definition and properties of
  {L}agrangian coherent structures from finite-time {L}yapunov exponents in
  two-dimensional aperiodic flows,'' \emph{Physica D: Nonlinear Phenomena},
  vol. 212, no. 3-4, pp. 271--304, 2005.

\bibitem{Green2007jfm}
M.~A. Green, C.~W. Rowley, and G.~Haller, ``Detection of {Lagrangian} coherent
  structures in {3D} turbulence.'' \emph{Journal of Fluid Mechanics}, vol. 572,
  pp. 111--120, 2007.

\bibitem{brunton2010chaos}
S.~L. Brunton and C.~W. Rowley, ``Fast computation of finite-time {L}yapunov
  exponent fields for unsteady flows,'' \emph{Chaos: An Interdisciplinary
  Journal of Nonlinear Science}, vol.~20, no.~1, p. 017503, 2010.

\bibitem{lipinski:2010}
D.~Lipinski and K.~Mohseni, ``A ridge tracking algorithm and error estimate for
  efficient computation of {L}agrangian coherent structures,'' \emph{Chaos},
  vol.~20, p. 017503, 2010.

\bibitem{haller2015arfm}
G.~Haller, ``Lagrangian coherent structures,'' \emph{Annual Review of Fluid
  Mechanics}, vol.~47, pp. 137--162, 2015.

\bibitem{haller2002pof}
------, ``Lagrangian coherent structures from approximate velocity data,''
  \emph{Physics of fluids}, vol.~14, no.~6, pp. 1851--1861, 2002.

\bibitem{haller:05}
------, ``An objective definition of a vortex,'' \emph{Journal of Fluid
  Mechanics}, vol. 525, pp. 1--26, 2005.

\bibitem{shadden2009correlation}
S.~C. Shadden, F.~Lekien, J.~D. Paduan, F.~P. Chavez, and J.~E. Marsden, ``The
  correlation between surface drifters and coherent structures based on
  high-frequency radar data in {M}onterey {B}ay,'' \emph{Deep Sea Research Part
  II: Topical Studies in Oceanography}, vol.~56, no. 3-5, pp. 161--172, 2009.

\bibitem{shadden2011lagrangian}
S.~C. Shadden, ``Lagrangian coherent structures,'' \emph{Transport and Mixing
  in Laminar Flows: From Microfluidics to Oceanic Currents}, pp. 59--89, 2011.

\bibitem{sudharsan2016lagrangian}
M.~Sudharsan, S.~L. Brunton, and J.~J. Riley, ``Lagrangian coherent structures
  and inertial particle dynamics,'' \emph{Physical Review E}, vol.~93, no.~3,
  p. 033108, 2016.

\bibitem{wilson2009lagrangian}
M.~M. Wilson, J.~Peng, J.~O. Dabiri, and J.~D. Eldredge, ``Lagrangian coherent
  structures in low {Reynolds} number swimming,'' \emph{Journal of Physics:
  Condensed Matter}, vol.~21, no.~20, p. 204105, 2009.

\bibitem{shadden2008characterization}
S.~C. Shadden and C.~A. Taylor, ``Characterization of coherent structures in
  the cardiovascular system,'' \emph{Annals of Biomedical Engineering},
  vol.~36, no.~7, pp. 1152--1162, 2008.

\bibitem{forgoston2011maximal}
E.~Forgoston, S.~Bianco, L.~B. Shaw, and I.~B. Schwartz, ``Maximal sensitive
  dependence and the optimal path to epidemic extinction,'' \emph{Bulletin of
  mathematical biology}, vol.~73, no.~3, pp. 495--514, 2011.

\bibitem{tallapragada2011lagrangian}
P.~Tallapragada, S.~D. Ross, and D.~G. Schmale~III, ``Lagrangian coherent
  structures are associated with fluctuations in airborne microbial
  populations,'' \emph{Chaos: An Interdisciplinary Journal of Nonlinear
  Science}, vol.~21, no.~3, p. 033122, 2011.

\bibitem{rockwood2019practical}
M.~P. Rockwood, T.~Loiselle, and M.~A. Green, ``Practical concerns of
  implementing a finite-time {L}yapunov exponent analysis with under-resolved
  data,'' \emph{Experiments in Fluids}, vol.~60, no.~4, p.~74, 2019.

\bibitem{rockwood2019real}
M.~P. Rockwood and M.~A. Green, ``Real-time identification of vortex shedding
  in the wake of a circular cylinder,'' \emph{AIAA Journal}, vol.~57, no.~1,
  pp. 223--238, 2019.

\bibitem{kularatne2016time}
D.~Kularatne, S.~Bhattacharya, and M.~A. Hsieh, ``Time and energy optimal path
  planning in general flows.'' in \emph{Robotics: Science and Systems}, 2016.

\bibitem{rao2009large}
D.~Rao and S.~B. Williams, ``Large-scale path planning for underwater gliders
  in ocean currents,'' in \emph{Australasian conference on robotics and
  automation (ACRA)}, 2009, pp. 2--4.

\bibitem{subramani2016energy}
D.~N. Subramani and P.~F. Lermusiaux, ``Energy-optimal path planning by
  stochastic dynamically orthogonal level-set optimization,'' \emph{Ocean
  Modelling}, vol. 100, pp. 57--77, 2016.

\bibitem{yilmaz2008path}
N.~K. Yilmaz, C.~Evangelinos, P.~F. Lermusiaux, and N.~M. Patrikalakis, ``Path
  planning of autonomous underwater vehicles for adaptive sampling using mixed
  integer linear programming,'' \emph{IEEE Journal of Oceanic Engineering},
  vol.~33, no.~4, pp. 522--537, 2008.

\bibitem{lermusiaux2007adaptive}
P.~F. Lermusiaux, ``Adaptive modeling, adaptive data assimilation and adaptive
  sampling,'' \emph{Physica D: Nonlinear Phenomena}, vol. 230, no. 1-2, pp.
  172--196, 2007.

\bibitem{bhatta2005coordination}
P.~Bhatta, E.~Fiorelli, F.~Lekien, N.~E. Leonard, D.~Paley, F.~Zhang,
  R.~Bachmayer, R.~E. Davis, D.~M. Fratantoni, and R.~Sepulchre, ``Coordination
  of an underwater glider fleet for adaptive ocean sampling,'' in \emph{Proc.
  International Workshop on Underwater Robotics, Int. Advanced Robotics
  Programmed (IARP), Genoa, Italy}, 2005.

\bibitem{leonard2007collective}
N.~E. Leonard, D.~A. Paley, F.~Lekien, R.~Sepulchre, D.~M. Fratantoni, and
  R.~E. Davis, ``Collective motion, sensor networks, and ocean sampling,''
  \emph{Proceedings of the IEEE}, vol.~95, no.~1, pp. 48--74, 2007.

\bibitem{fiorelli2006multi}
E.~Fiorelli, N.~E. Leonard, P.~Bhatta, D.~A. Paley, R.~Bachmayer, and D.~M.
  Fratantoni, ``Multi-{AUV} control and adaptive sampling in {M}onterey
  {B}ay,'' \emph{IEEE Journal of Oceanic Engineering}, vol.~31, no.~4, pp.
  935--948, 2006.

\bibitem{leonard2001model}
N.~E. Leonard and J.~G. Graver, ``Model-based feedback control of autonomous
  underwater gliders,'' \emph{IEEE Journal of Oceanic Engineering}, vol.~26,
  no.~4, pp. 633--645, 2001.

\bibitem{lipinski2010cooperative}
D.~Lipinski and K.~Mohseni, ``Cooperative control of a team of unmanned
  vehicles using smoothed particle hydrodynamics,'' in \emph{AIAA Guidance,
  Navigation, and Control Conference}, 2010, p. 8316.

\bibitem{lipinski2014feasible}
------, ``Feasible area coverage of a hurricane using micro-aerial vehicles,''
  in \emph{AIAA Atmospheric Flight Mechanics Conference}, 2014, p. 0894.

\bibitem{lipinski2011master}
------, ``A master-slave fluid cooperative control algorithm for optimal
  trajectory planning,'' in \emph{2011 IEEE International Conference on
  Robotics and Automation}.\hskip 1em plus 0.5em minus 0.4em\relax IEEE, 2011,
  pp. 3347--3351.

\bibitem{krishna2022finite}
K.~Krishna, Z.~Song, and S.~L. Brunton, ``Finite-horizon, energy-efficient
  trajectories in unsteady flows,'' \emph{Proceedings of the Royal Society A},
  vol. 478, no. 2258, p. 20210255, 2022.

\bibitem{gunnarson2021learning}
P.~Gunnarson, I.~Mandralis, G.~Novati, P.~Koumoutsakos, and J.~O. Dabiri,
  ``Learning efficient navigation in vortical flow fields,'' \emph{arXiv
  preprint arXiv:2102.10536}, 2021.

\bibitem{bifarale2019}
\BIBentryALTinterwordspacing
L.~Biferale, F.~Bonaccorso, M.~Buzzicotti, P.~Clark Di~Leoni, and
  K.~Gustavsson, ``Zermelo’s problem: Optimal point-to-point navigation in
  2{D} turbulent flows using reinforcement learning,'' \emph{Chaos: An
  Interdisciplinary Journal of Nonlinear Science}, vol.~29, no.~10, p. 103138,
  2019. [Online]. Available: \url{https://doi.org/10.1063/1.5120370}
\BIBentrySTDinterwordspacing

\bibitem{buzzicotti2021optimal}
M.~Buzzicotti, L.~Biferale, F.~Bonaccorso, P.~C. di~Leoni, and K.~Gustavsson,
  ``Optimal control of point-to-point navigation in turbulent time-dependent
  flows using reinforcement learning,'' 2021.

\bibitem{jiao2021learning}
Y.~Jiao, F.~Ling, S.~Heydari, N.~Heess, J.~Merel, and E.~Kanso, ``Learning to
  swim in potential flow,'' \emph{Physical Review Fluids}, vol.~6, no.~5, p.
  050505, 2021.

\bibitem{ao2023}
\BIBentryALTinterwordspacing
A.~Xu, H.-L. Wu, and H.-D. Xi, ``Long-distance migration with minimal energy
  consumption in a thermal turbulent environment,'' \emph{Phys. Rev. Fluids},
  vol.~8, p. 023502, Feb 2023. [Online]. Available:
  \url{https://link.aps.org/doi/10.1103/PhysRevFluids.8.023502}
\BIBentrySTDinterwordspacing

\bibitem{inanc2005optimal}
T.~Inanc, S.~C. Shadden, and J.~E. Marsden, ``Optimal trajectory generation in
  ocean flows,'' in \emph{Proceedings of the 2005, American Control Conference,
  2005.}\hskip 1em plus 0.5em minus 0.4em\relax IEEE, 2005, pp. 674--679.

\bibitem{zhang2008optimal}
W.~Zhang, T.~Inanc, S.~Ober-Blobaum, and J.~E. Marsden, ``Optimal trajectory
  generation for a glider in time-varying {2D} ocean flows {B}-spline model,''
  in \emph{2008 IEEE International Conference on Robotics and
  Automation}.\hskip 1em plus 0.5em minus 0.4em\relax IEEE, 2008, pp.
  1083--1088.

\bibitem{senatore2008fuel}
C.~Senatore and S.~D. Ross, ``Fuel-efficient navigation in complex flows,'' in
  \emph{2008 American Control Conference}.\hskip 1em plus 0.5em minus
  0.4em\relax IEEE, 2008, pp. 1244--1248.

\bibitem{heckman2016controlling}
C.~R. Heckman, M.~A. Hsieh, and I.~B. Schwartz, ``Controlling basin breakout
  for robots operating in uncertain flow environments,'' in \emph{Experimental
  Robotics}.\hskip 1em plus 0.5em minus 0.4em\relax Springer, 2016, pp.
  561--576.

\bibitem{garcia1989model}
C.~E. Garcia, D.~M. Prett, and M.~Morari, ``Model predictive control: {T}heory
  and practice---{A} survey,'' \emph{Automatica}, vol.~25, no.~3, pp. 335--348,
  1989.

\bibitem{camacho2013model}
E.~F. Camacho and C.~B. Alba, \emph{Model Predictive Control}.\hskip 1em plus
  0.5em minus 0.4em\relax Springer Science \& Business Media, 2013.

\bibitem{sutton2018reinforcement}
R.~S. Sutton and A.~G. Barto, \emph{Reinforcement learning: An introduction},
  2018.

\bibitem{recht2019tour}
B.~Recht, ``A tour of reinforcement learning: The view from continuous
  control,'' \emph{Annual Review of Control, Robotics, and Autonomous Systems},
  vol.~2, pp. 253--279, 2019.

\bibitem{Brunton2022book}
S.~L. Brunton and J.~N. Kutz, \emph{Data-Driven Science and Engineering:
  Machine Learning, Dynamical Systems, and Control}, 2nd~ed.\hskip 1em plus
  0.5em minus 0.4em\relax Cambridge University Press, 2022.

\bibitem{peng2009transport}
J.~Peng and J.~Dabiri, ``Transport of inertial particles by lagrangian coherent
  structures: application to predator--prey interaction in jellyfish feeding,''
  \emph{Journal of Fluid Mechanics}, vol. 623, pp. 75--84, 2009.

\bibitem{kelley2013lagrangian}
D.~H. Kelley, M.~R. Allshouse, and N.~T. Ouellette, ``Lagrangian coherent
  structures separate dynamically distinct regions in fluid flows,''
  \emph{Physical Review E}, vol.~88, no.~1, p. 013017, 2013.

\bibitem{olascoaga2006persistent}
M.~J. Olascoaga, I.~Rypina, M.~G. Brown, F.~J. Beron-Vera, H.~Ko{\c{c}}ak,
  L.~E. Brand, G.~Halliwell, and L.~K. Shay, ``Persistent transport barrier on
  the {West Florida Shelf},'' \emph{Geophysical research letters}, vol.~33,
  no.~22, 2006.

\bibitem{beron2008oceanic}
F.~J. Beron-Vera, M.~J. Olascoaga, and G.~Goni, ``Oceanic mesoscale eddies as
  revealed by {L}agrangian coherent structures,'' \emph{Geophysical Research
  Letters}, vol.~35, no.~12, 2008.

\bibitem{beron2015dissipative}
F.~J. Beron-Vera, M.~J. Olascoaga, G.~Haller, M.~Farazmand, J.~Tri{\~n}anes,
  and Y.~Wang, ``Dissipative inertial transport patterns near coherent
  {L}agrangian eddies in the ocean,'' \emph{Chaos: An Interdisciplinary Journal
  of Nonlinear Science}, vol.~25, no.~8, p. 087412, 2015.

\bibitem{Lekien2005physicad}
F.~Lekien, C.~Coulliette, A.~J. Mariano, E.~H. Ryan, L.~K. Shay, G.~Haller, and
  J.~E. Marsden, ``Pollution release tied to invariant manifolds: {A} case
  study for the coast of {F}lorida,'' \emph{Physica D}, vol. 210, pp. 1--20,
  2005.

\bibitem{franco:2007}
E.~Franco, D.~N. Pekarek, J.~Peng, and J.~O. Dabiri, ``Geometry of unsteady
  fluid transport during fluid-structure interactions,'' \emph{Journal of Fluid
  Mechanics}, vol. 589, pp. 125--145, 2007.

\bibitem{Padberg2007njp}
K.~Padberg, T.~Hauff, F.~Jenko, and O.~Junge, ``{Lagrangian} structures and
  transport in turbulent magnetized plasmas,'' \emph{New Journal of Physics},
  vol.~9, p. 400, 2007.

\bibitem{Mathur2007prl}
M.~Mathur, G.~Haller, T.~Peacock, J.~E. Ruppert-Felsot, and H.~L. Swinney,
  ``Uncovering the {Lagrangian} skeleton of turbulence,'' \emph{Physical Review
  Letters}, vol.~98, pp. 144\,502--1--144\,502--4, 2007.

\bibitem{Peng2008jeb}
J.~Peng and J.~O. Dabiri, ``The `upstream wake' of swimming and flying animals
  and its correlation with propulsive efficiency,'' \emph{The Journal of
  Experimental Biology}, vol. 211, pp. 2669--2677, 2008.

\bibitem{rockwood2016detecting}
M.~P. Rockwood, K.~Taira, and M.~A. Green, ``Detecting vortex formation and
  shedding in cylinder wakes using {L}agrangian coherent structures,''
  \emph{AIAA Journal}, vol.~55, no.~1, pp. 15--23, 2016.

\bibitem{stengel1994optimal}
R.~F. Stengel, \emph{Optimal control and estimation}.\hskip 1em plus 0.5em
  minus 0.4em\relax Courier Corporation, 1994.

\bibitem{todorov2006optimal}
E.~Todorov \emph{et~al.}, ``Optimal control theory,'' \emph{Bayesian brain:
  probabilistic approaches to neural coding}, pp. 268--298, 2006.

\bibitem{JiangS:95a}
S.~Jiang, F.-f. Jin, and M.~Ghil, ``Multiple equilibria, periodic, and
  aperiodic solutions in a wind-driven, double-gyre, shallow-water model,''
  \emph{Journal of Physical Oceanography}, vol.~25, no.~5, pp. 764--786, 1995.

\bibitem{SpeichS:94a}
S.~Speich and M.~Ghil, ``Interannual variability of the mid-latitude oceans: A
  new source of climate variability,'' \emph{Sistema Terra}, vol.~3, no.~3, p.
  459, 1994.

\bibitem{SpeichS:95a}
S.~Speich, H.~Dijkstra, and M.~Ghil, ``Successive bifurcations in a
  shallow-water model applied to the wind-driven ocean circulation,''
  \emph{Nonlinear Processes in Geophysics}, vol.~2, no. 3/4, pp. 241--268,
  1995.

\bibitem{Andersson2019}
J.~A.~E. Andersson, J.~Gillis, G.~Horn, J.~B. Rawlings, and M.~Diehl,
  ``{CasADi} -- {A} software framework for nonlinear optimization and optimal
  control,'' \emph{Mathematical Programming Computation}, vol.~11, no.~1, pp.
  1--36, 2019.

\bibitem{rawlings2015}
M.~J. Risbeck and J.~B. Rawlings, ``{MPCTools}: Nonlinear model predictive
  control tools for {CasADi} ({Python} interface),'' 2015,
  {https://bitbucket.org/rawlings-group/mpc-tools-casadi}.

\bibitem{stable-baselines3}
\BIBentryALTinterwordspacing
A.~Raffin, A.~Hill, A.~Gleave, A.~Kanervisto, M.~Ernestus, and N.~Dormann,
  ``Stable-baselines3: Reliable reinforcement learning implementations,''
  \emph{Journal of Machine Learning Research}, vol.~22, no. 268, pp. 1--8,
  2021. [Online]. Available: \url{http://jmlr.org/papers/v22/20-1364.html}
\BIBentrySTDinterwordspacing

\bibitem{Talley:19a}
\BIBentryALTinterwordspacing
L.~D. Talley, I.~Rosso, I.~Kamenkovich, M.~R. Mazloff, J.~Wang, E.~Boss, A.~R.
  Gray, K.~S. Johnson, R.~M. Key, S.~C. Riser, N.~L. Williams, and J.~L.
  Sarmiento, ``{Southern Ocean} biogeochemical float deployment strategy, with
  example from the {Greenwich Meridian} line ({GO-SHIP A12}),'' \emph{Journal
  of Geophysical Research: Oceans}, vol. 124, no.~1, pp. 403--431, 2019.
  [Online]. Available:
  \url{https://agupubs.onlinelibrary.wiley.com/doi/abs/10.1029/2018JC014059}
\BIBentrySTDinterwordspacing

\bibitem{wiggins2013chaotic}
S.~Wiggins, \emph{Chaotic transport in dynamical systems}.\hskip 1em plus 0.5em
  minus 0.4em\relax Springer Science \& Business Media, 2013, vol.~2.

\bibitem{rom1990analytical}
V.~Rom-Kedar, A.~Leonard, and S.~Wiggins, ``An analytical study of transport,
  mixing and chaos in an unsteady vortical flow,'' \emph{Journal of Fluid
  Mechanics}, vol. 214, pp. 347--394, 1990.

\bibitem{sastry2013nonlinear}
S.~Sastry, \emph{Nonlinear systems: analysis, stability, and control}.\hskip
  1em plus 0.5em minus 0.4em\relax Springer Science \& Business Media, 2013,
  vol.~10.

\bibitem{lillicrap2015continuous}
T.~P. Lillicrap, J.~J. Hunt, A.~Pritzel, N.~Heess, T.~Erez, Y.~Tassa,
  D.~Silver, and D.~Wierstra, ``Continuous control with deep reinforcement
  learning,'' \emph{arXiv preprint arXiv:1509.02971}, 2015.

\bibitem{rakovic2007optimized}
S.~V. Rakovi{\'c}, E.~C. Kerrigan, D.~Q. Mayne, and K.~I. Kouramas, ``Optimized
  robust control invariance for linear discrete-time systems: Theoretical
  foundations,'' \emph{Automatica}, vol.~43, no.~5, pp. 831--841, 2007.

\bibitem{wiggert2022navigating}
M.~Wiggert, M.~Doshi, P.~F. Lermusiaux, and C.~J. Tomlin, ``Navigating
  underactuated agents by hitchhiking forecast flows,'' in \emph{2022 IEEE 61st
  Conference on Decision and Control (CDC)}.\hskip 1em plus 0.5em minus
  0.4em\relax IEEE, 2022, pp. 2417--2424.

\bibitem{doshi2023energy}
M.~M. Doshi, M.~S. Bhabra, and P.~F. Lermusiaux, ``Energy--time optimal path
  planning in dynamic flows: Theory and schemes,'' \emph{Computer Methods in
  Applied Mechanics and Engineering}, vol. 405, p. 115865, 2023.

\bibitem{guckenheimer2004fast}
J.~Guckenheimer and A.~Vladimirsky, ``A fast method for approximating invariant
  manifolds,'' \emph{SIAM Journal on Applied Dynamical Systems}, vol.~3, no.~3,
  pp. 232--260, 2004.

\bibitem{guckenheimer2015invariant}
J.~Guckenheimer, B.~Krauskopf, H.~M. Osinga, and B.~Sandstede, ``Invariant
  manifolds and global bifurcations,'' \emph{Chaos: An Interdisciplinary
  Journal of Nonlinear Science}, vol.~25, no.~9, p. 097604, 2015.

\bibitem{gonzalez2016finite}
D.~Gonz{\'a}lez, R.~Speth, D.~Gaitonde, and M.~Lewis, ``Finite-time lyapunov
  exponent-based analysis for compressible flows,'' \emph{Chaos: An
  Interdisciplinary Journal of Nonlinear Science}, vol.~26, no.~8, p. 083112,
  2016.

\end{thebibliography}
 }

\end{document}